\documentclass[12pt,reqno]{amsart}
\usepackage{amssymb}
\usepackage{amscd}
\usepackage{bm}

\newtheorem{theorem}{Theorem}[section]
\newtheorem{corollary}[theorem]{Corollary}
\newtheorem{lemma}[theorem]{Lemma}
\newtheorem{proposition}[theorem]{Proposition}

\newtheorem{definition}[theorem]{Definition}
\newtheorem{remark}[theorem]{Remark}

\newtheorem{example}[theorem]{Example}
\numberwithin{equation}{section}
\begin{document}
\title{Witten deformation and the equivariant index}
\author{Igor Prokhorenkov}
\author{Ken Richardson}
\address{Department of Mathematics\\
Texas Christian University\\
Box 298900 \\
Fort Worth, Texas 76129}
\email{i.prokhorenkov@tcu.edu\\
k.richardson@tcu.edu}
\date{October, 2006}
\maketitle

\begin{abstract}
Let $M$ be a compact Riemannian manifold endowed with an isometric action of
a compact, connected Lie group. The method of the Witten deformation is used
to compute the virtual representation-valued equivariant index of a
transversally elliptic, first order differential operator on $M$. The
multiplicities of irreducible representations in the index are expressed in
terms of local quantities associated to the isolated singular points of an
equivariant bundle map that is locally Clifford multiplication by a Killing
vector field near these points.
\end{abstract}

\section{Introduction}

The purpose of this paper is to compute the equivariant index multiplicities
of an equivariant, transversally elliptic operator on a compact $G$%
-manifold, where $G$ is a compact, connected Lie group. We use the method of
Witten deformation to express the index in terms of combinatorial data
associated to a given equivariant bundle map.

We start by establishing notation and reviewing the definitions of various
types of equivariant indices associated to first order, transversally
elliptic differential operators. In Section \ref{content} we explain our
application of the Witten deformation technique for calculating these
equivariant indices and discuss the main results of this paper.

\subsection{Types of equivariant indices}

Suppose that a compact Lie group $G$ acts by isometries on a compact,
connected Riemannian manifold $M$, and let $E=E^{+}\oplus E^{-}$ be a
graded, $G$-equivariant, Hermitian vector bundle over $M$. We consider a
first order $G$-equivariant differential operator $D^{+}:$ $\Gamma \left(
M,E^{+}\right) \rightarrow \Gamma \left( M,E^{-}\right) $ which is elliptic
merely in the directions transversal to the orbits of $G$, and let $D^{-}$
be the formal adjoint of $D^{+}$. Then the operator $D^{+}$ belongs to the
class of \emph{transversally elliptic differential operators} introduced by
M. Atiyah in \cite{A}. In this paper, we will assume for the most part that $%
G$ is connected and that the operator $D^{+}$ is in addition transversally
elliptic with respect to the action of a maximal torus in $G$.

The group $G$ acts in a natural way on $\Gamma \left( M,E^{\pm }\right) $,
and the (possibly infinite-dimensional) subspaces $\ker \left( D^{+}\right) $
and $\ker \left( D^{-}\right) $ are $G$-invariant subspaces. Thus, each of $%
\Gamma \left( M,E^{\pm }\right) $, $\ker \left( D^{+}\right) $, and $\ker
\left( D^{-}\right) $ decomposes as a direct sum of irreducible
representation spaces. Let $\rho :G\rightarrow \mathrm{End}\left( V_{\rho
}\right) $ be an irreducible unitary representation of $G$, and let $\chi
_{\rho }:G\rightarrow \mathbb{C}$ be its character; that is, $\chi _{\rho
}\left( g\right) =\mathrm{tr}\left( \rho \left( g\right) \right) $. By the
Peter-Weyl Theorem, the functions $\left\{ \chi _{\rho }\right\} _{\rho }$
are eigenfunctions of the Laplacian on $G$ and form an orthonormal set in $%
L^{2}\left( G\right) $ with the normalized, biinvariant metric. Let $\Gamma
\left( M,E^{\pm }\right) ^{\rho }$ be the subspace of sections that is the
direct sum of the irreducible $G$-representation subspaces of $\Gamma \left(
M,E^{\pm }\right) $ corresponding to representations that are unitarily
equivalent to $\rho $. It can be shown that the operator

\begin{equation*}
D^{+}:\Gamma \left( M,E^{+}\right) ^{\rho }\rightarrow \Gamma \left(
M,E^{-}\right) ^{\rho }
\end{equation*}%
can be extended to a Fredholm operator between the appropriate Sobolev
spaces, so that each irreducible representation of $G$ appears with finite
multiplicity in $\ker D^{\pm }$. Let $a_{\rho }^{\pm }\in \mathbb{Z}^{+}$ be
the multiplicity of $\rho $ in $\ker \left( D^{\pm }\right) $.

As in \cite{A}, we define the virtual representation-valued index of $D$ to
be%
\begin{equation*}
\mathrm{ind}^{G}\left( D^{+}\right) :=\sum_{\rho }\left( a_{\rho
}^{+}-a_{\rho }^{-}\right) \left[ \rho \right] ,
\end{equation*}%
where $\left[ \rho \right] $ denotes the equivalence class of the
irreducible representation $\rho $. The index multiplicity is 
\begin{equation*}
\mathrm{ind}^{\rho }\left( D^{+}\right) :=a_{\rho }^{+}-a_{\rho }^{-}=\frac{1%
}{\dim V_{\rho }}\mathrm{ind}\left( \left. D^{+}\right\vert _{\Gamma \left(
M,E^{+}\right) ^{\rho }\rightarrow \Gamma \left( M,E^{-}\right) ^{\rho
}}\right) .
\end{equation*}%
In particular, if $\rho _{0}$ is the trivial representation of $G$, then 
\begin{equation*}
\mathrm{ind}^{\rho _{0}}\left( D^{+}\right) =\mathrm{ind}\left( \left.
D^{+}\right\vert _{\Gamma \left( M,E^{+}\right) ^{G}\rightarrow \Gamma
\left( M,E^{-}\right) ^{G}}\right) ,
\end{equation*}%
where the superscript $G$ implies restriction to $G$-invariant sections.

The relationship between the index multiplicities and Atiyah's equivariant
distribution-valued index $\mathrm{ind}_{g}\left( D^{+}\right) $ is as
follows. The virtual character $\mathrm{ind}_{g}\left( D^{+}\right) $ is
given by (see \cite{A}) 
\begin{eqnarray*}
\mathrm{ind}_{g}\left( D^{+}\right) &:&=\text{\textquotedblleft }\mathrm{tr}%
\left( \left. g\right\vert _{\ker D^{+}}\right) -\mathrm{tr}\left( \left.
g\right\vert _{\ker D^{-}}\right) \text{\textquotedblright } \\
&=&\sum_{\rho }\mathrm{ind}^{\rho }\left( D^{+}\right) \chi _{\rho }\left(
g\right) \in \mathcal{D}\left( G\right) ,
\end{eqnarray*}%
where $\mathcal{D}\left( G\right) $ is the set of distributions on $G$.
Since $\ker D^{+}$ and $\ker D^{-}$ are in general infinite-dimensional, the
sum above does not always converge, but it makes sense as a distribution on $%
G$. That is, if $dg$ is the normalized, biinvariant volume form on $G$, and
if $\phi =\sum c_{\rho }\chi _{\rho }\in C^{\infty }\left( G\right) $, then 
\begin{eqnarray*}
\mathrm{ind}_{\left( \bullet \right) }\left( D^{+}\right) \left( \phi
\right) &=&\text{\textquotedblleft }\int_{G}\phi \left( g\right) ~\overline{%
\mathrm{ind}_{g}\left( D^{+}\right) }~dg\text{\textquotedblright } \\
&=&\sum_{\rho }\mathrm{ind}^{\rho }\left( D^{+}\right) \int \phi \left(
g\right) ~\overline{\chi _{\rho }\left( g\right) }~dg=\sum_{\rho }\mathrm{ind%
}^{\rho }\left( D^{+}\right) c_{\rho },
\end{eqnarray*}%
an expression which converges because the coefficients $c_{\rho }$ are
rapidly decreasing and $\mathrm{ind}^{\rho }\left( D^{+}\right) $ grows at
most polynomially as $\rho $ varies over the irreducible representations of $%
G$. From this calculation, we see that the multiplicities determine Atiyah's
distributional index. Conversely, let $\alpha :G\rightarrow \mathrm{End}%
\left( V_{\alpha }\right) $ be an irreducible unitary representation. Then 
\begin{equation*}
\mathrm{ind}_{\left( \bullet \right) }\left( D^{+}\right) \left( \chi
_{\alpha }\right) =\sum_{\rho }\mathrm{ind}^{\rho }\left( D^{+}\right) \int
\chi _{\alpha }\left( g\right) \overline{\chi _{\rho }\left( g\right) }\,dg=%
\mathrm{ind}^{\alpha }D^{+},
\end{equation*}%
so that complete knowledge of the equivariant distributional index is
equivalent to knowing all of the multiplicities $\mathrm{ind}^{\rho }\left(
D^{+}\right) $. Because the operator $\left. D^{+}\right\vert _{\Gamma
\left( M,E^{+}\right) ^{\rho }\rightarrow \Gamma \left( M,E^{-}\right)
^{\rho }}$ is Fredholm, all of the indices $\mathrm{ind}^{G}\left(
D^{+}\right) $ , $\mathrm{ind}_{g}\left( D^{+}\right) $, and $\mathrm{ind}%
^{\rho }\left( D^{+}\right) $ depend only on the equivariant homotopy class
of the principal transverse symbol of $D^{+}$.

\subsection{Content of the paper: applications of Witten deformation to
equivariant index theory\label{content}}

About 25 years ago E. Witten \cite{Wit1} introduced a new way of proving
Morse inequalities based on a deformation of the de Rham complex. His ideas
were fruitfully applied in many specific situations. The purpose of this
paper is to utilize this method to prove an explicit formula for the index $%
\mathrm{ind}^{\rho }\left( D^{+}\right) $ in terms of data associated to the
singular set of an equivariant bundle map $Z:E\rightarrow E$. In this paper,
we require that the singularities, if they exist, are isolated and that the
map $Z$ has the form%
\begin{equation*}
Z=c\left( iV\right)
\end{equation*}%
near each singular point, where $c$ denotes a locally defined Clifford
multiplication and $V$ is a Killing vector field. Witten used a similar
approach in \cite{Wit3} to prove the Atiyah-Hirzebruch vanishing theorem
(see \cite{AH}) by showing $\mathrm{ind}^{\rho }\left( D^{+}\right) =0$ if $%
D^{+}$ is the Dirac operator on spinors and $G=S^{1}$; in his argument $%
Z=c\left( i\partial _{\theta }\right) $ globally. It should be mentioned
that this idea is related to Atiyah's earlier method of \textquotedblleft
pushing a symbol,\textquotedblright\ to extract information about the
distribution-valued equivariant index near fixed points of a torus action
(see \cite[Chapter 6]{A}).

Let $D=\left( 
\begin{array}{cc}
0 & D^{-} \\ 
D^{+} & 0%
\end{array}%
\right) $. We consider the following family of transversally elliptic
operators, depending on a real parameter $s$ :%
\begin{eqnarray*}
D_{s} &=&D+sZ,\text{ so that} \\
D_{s}^{2} &=&D^{2}+s\left( DZ+ZD\right) +s^{2}Z^{2}
\end{eqnarray*}%
We want to study the spectral asymptotics of this family as $s\rightarrow
\infty $. Unlike most other applications of Witten deformation where the
operator $B=DZ+ZD$ is bounded (see \cite{Wit1}, \cite{Pro-Rich}), in this
paper the operator $B$ is first order (at least near singular points) and
thus unbounded. In order to circumvent this difficulty, we require that the
restriction of the $B$ to $\Gamma \left( M,E\right) ^{\rho }$ is a bundle
map, which is indeed true if the first order part of $B$ is a tangential
derivative. In Section \ref{localizationSection}, we extend the localization
theorem of Shubin (\cite{Shu1}) to the setting of transversally elliptic
operators. This result allows us to reduce the computation of $\mathrm{ind}%
^{\rho }\left( D^{+}\right) $ to investigating the spectrum of a certain
model operator at each singular point of $Z$.

In Section \ref{EquivariantWitten}, we restrict to the case where $G$ is a
torus, and we compute the spectral asymptotics of the operator $\frac{1}{s}%
D_{s}^{2}$ as $s\rightarrow \infty $ in terms of local information at each
singular point of $Z$. The main result of the section is Theorem \ref%
{ModelSpectrumTheorem}.

In Section \ref{EquIndexThmSection}, we apply Theorem \ref%
{ModelSpectrumTheorem} to evaluate the index $\mathrm{ind}^{\rho }\left(
D^{+}\right) $ in the case where $G$ is a torus. The main result of the
paper is the formula for this index in the Transverse Index Theorem, Theorem %
\ref{transverseGIndexTheorem}. In Section \ref{torusindices}, we show that
for any compact, connected Lie group $G$, the index $\mathrm{ind}^{\rho
}\left( D^{+}\right) $ can be expressed in terms of the corresponding
indices for its maximal torus, as long as the relevant torus multiplicities
are finite (as in the case where $D^{+}$ is also transversally elliptic with
respect to the torus action).

Finally, in Section \ref{ExampleSection}, we demonstrate applications of
Theorem \ref{transverseGIndexTheorem} to the signature and de Rham operators
on $G$-manifolds and to a specific transversally elliptic operator on the
sphere. These investigations yield an interesting new identity involving
Killing vector fields on $G$-manifolds along with new proofs of other known
identities; see Proposition \ref{killingProp}. We also apply the theory in
Section \ref{torusindices} to an example of an $SU\left( 2\right) $-action
on a sphere.

\subsection{Historical Comments}

A large body of work over the last twenty years has yielded theorems that
express $\mathrm{ind}_{g}\left( D^{+}\right) $ and the corresponding local
heat kernel supertrace in terms of topological and geometric quantities (as
in the Atiyah-Segal-Singer index theorem for elliptic operators or the
Berline-Vergne Theorem for transversally elliptic operators --- see \cite%
{ASe},\cite{Be-V1},\cite{Be-V2}). The problem of expressing $\mathrm{ind}%
^{\rho }\left( D^{+}\right) $ explicitly as a sum of topological or
geometric quantities which are determined at the different strata of the $G$%
-manifold $M$ is addressed in the paper \cite{BruKRich}. The special cases
where $G$ is finite or when all of the isotropy groups have the same
dimension were solved by M. Atiyah in \cite{A}, and it turns out both of
these are special cases of the Orbifold Index Theorem by T. Kawasaki (see 
\cite{Kawas2}). In the case when $D^{+}$ is elliptic, the Atiyah-Bott fixed
point formula may be used to calculate the equivariant indices corresponding
to a torus action from fixed point data, as in this paper (see \cite{AB1},%
\cite{AB2}). Much work has also been done on symplectic manifolds, where the
local data comes from the critical set of the moment map. For example, see 
\cite{TZ} for an analytic proof of the Guillemin-Sternberg conjecture (\cite%
{G-S}). Also, see \cite{Par} and \cite{Bra} for another Witten deformation
approach to finding the equivariant index of a specific transversally
elliptic symbol on a noncompact manifold.

\section{Equivariant Localization\label{localizationSection}}

Suppose a compact Lie group $G$ acts by isometries on a closed, oriented
Riemannian manifold $M$ of dimension $2n$. Let $E$ be a $G$-equivariant
Hermitian bundle over $M$. Let $\rho $ be an irreducible representation of $%
G $, and let $\Gamma \left( M,E\right) ^{\rho }$ denote the space of
sections of $E$ of type $\rho $. For $s>0$, let $H_{s}:\Gamma \left(
M,E\right) \rightarrow \Gamma \left( M,E\right) $ be a transversally
elliptic, $G$-equivariant, essentially self-adjoint, second order
differential operator of the form 
\begin{equation*}
H_{s}=\frac{1}{s}A+B+sC,
\end{equation*}%
where

\begin{enumerate}
\item $A$ is a second order, transversally elliptic differential operator
with positive definite principal transverse symbol.

\item For each irreducible representation $\rho $ of $G$, $\left. B\right|
_{\Gamma \left( M,E\right) ^{\rho }}$ is a bundle map.

\item $C$ is a bundle map such that $C\left( x\right) \geq 0$ for all $x\in
M $, and at each point $\overline{x}$ where $C\left( \overline{x}\right) $
is singular, there exists $c>0$ such that 
\begin{eqnarray*}
C\left( \overline{x}\right) &=&0,\,\overline{x}g=\overline{x}\text{ for all }%
g\in G,\,\,\text{and} \\
C\left( x\right) &\geq &c\cdot d\left( x,\overline{x}\right) ^{2}\mathbf{1}
\end{eqnarray*}%
in a neighborhood of $\overline{x}$, where $d\left( x,\overline{x}\right) $
is the distance from $x$ to $\overline{x}$.

\item $A$ is elliptic in a neighborhood of each singular point of $C$.
\end{enumerate}

Let $H_{s}^{\rho }$ denote the restriction of $H_{s}$ to $\Gamma \left(
M,E\right) ^{\rho }$. For each $\rho $, the operator $H_{s}^{\rho }$ has
discrete spectrum (see \cite[p. 12-13]{A}); this implies that the spectrum
of $H_{s}$ consists of a discrete set of eigenvalues, although some
eigenvalues may have infinite multiplicities.

Near each singular point $\overline{x}$ of $C$, we choose coordinates $x=$ $%
\left( x_{1},...,x_{2n}\right) $ such that $\overline{x}$ corresponds to the
origin, $T_{\overline{x}}M=\mathbb{R}^{2n}$, and the volume form at the
origin is $dx_{1}...dx_{2n}$. We choose a trivialization of $E$ near $%
\overline{x}$ so that $A$, $B$, and $C$ become differential operators with
matrix coefficients. We define the model operator $K_{\overline{x}}^{\rho
}:\Gamma \left( \mathbb{R}^{2n},E_{\overline{x}}\right) ^{\rho }\rightarrow
\Gamma \left( \mathbb{R}^{2n},E_{\overline{x}}\right) ^{\rho }$ by 
\begin{eqnarray*}
K_{\overline{x}}^{\rho } &=&\widetilde{A}+\widetilde{B}^{\rho }+\widetilde{C}%
,\text{ where} \\
\widetilde{A} &=&\text{the principal part of }A\text{ at }\overline{x} \\
\widetilde{B}^{\rho } &=&\left. B\right\vert _{\Gamma \left( M,E\right)
^{\rho }}\left( \overline{x}\right) \\
\widetilde{C} &=&\sum x_{i}x_{j}\left( \nabla _{i}\nabla _{j}C\right) _{%
\overline{x}}=\text{the quadratic part of }C\text{ at }\overline{x},
\end{eqnarray*}%
where $\nabla $ is the induced connection on $E\otimes E^{\ast }$. It is
easy to check that $\widetilde{C}$ is independent of the coordinates chosen.
Let $dg$ denote differential of the action of $g\in G$ at $\overline{x}$, so
we write $dg:\mathbb{R}^{2n}\rightarrow \mathbb{R}^{2n}$. Let the action of $%
g$ on $\mathbb{R}^{2n}\times E_{\overline{x}}$ be defined as 
\begin{equation*}
\left( x,v_{\overline{x}}\right) g=\left( dg\left( x\right) ,g\cdot v_{%
\overline{x}}\right) .
\end{equation*}

\begin{lemma}
The operator $K_{\overline{x}}^{\rho }$ is equivariant with respect to this $%
G$-action.
\end{lemma}

\begin{proof}
Since $H_{s}$ is equivariant for each $s>0$, it is elementary to show that
each of the operators $A$, $B$, and $C$ is equivariant. Then the principal
symbol of $A$ is equivariant, and in particular the principal symbol of $A$
at $\overline{x}$ is $G$-invariant. Thus, $\widetilde{A}$ is $G$-invariant.
Next, since $\left. B\right| _{\Gamma \left( M,E\right) ^{\rho }}$ is
equivariant, its restriction $\widetilde{B}^{\rho }$ to $\overline{x}$ is
also. Finally, since $C$ is equivariant and the connection is equivariant,
it follows that $\widetilde{C}$ is equivariant.
\end{proof}

\begin{lemma}
For each irreducible representation $\rho $ of $G$ and each fixed point $%
\overline{x}$ of $G$, the operator $K_{\overline{x}}^{\rho }:\Gamma \left( 
\mathbb{R}^{2n},E_{\overline{x}}\right) ^{\rho }\rightarrow \Gamma \left( 
\mathbb{R}^{2n},E_{\overline{x}}\right) ^{\rho }$ has discrete spectrum.
\end{lemma}

\begin{proof}
Consider the extended operator $K_{\overline{x}}^{\rho }:\Gamma \left( 
\mathbb{R}^{2n},E_{\overline{x}}\right) \rightarrow \Gamma \left( \mathbb{R}%
^{2n},E_{\overline{x}}\right) $ . This operator is elliptic and essentially
self-adjoint, and the operator is bounded below by $\left( C_{1}+C_{2}\cdot
\left| x\right| ^{2}\right) \mathbf{1}$, where $C_{1}\in \mathbb{R}$ and $%
C_{2}>0 $. Since this bound goes to infinity as $x\rightarrow \pm \infty $,
the operator $K_{\overline{x}}^{\rho }-\left( C_{1}-1\right) \mathbf{1}$ has
a compact resolvent. Thus, the restriction of $K_{\overline{x}}^{\rho }$ to $%
\Gamma \left( \mathbb{R}^{2n},E_{\overline{x}}\right) ^{\rho }$ also has a
compact resolvent.
\end{proof}

We define the model operator $K^{\rho }$ by 
\begin{equation*}
K^{\rho }=\bigoplus_{\text{fixed point }\overline{x}}K_{\overline{x}}^{\rho
}.
\end{equation*}
Clearly, this operator has discrete spectrum. Let 
\begin{equation*}
\mu _{1}^{\rho }<\mu _{2}^{\rho }<\mu _{3}^{\rho }<...
\end{equation*}
be the distinct eigenvalues of $K^{\rho }$ with corresponding multiplicities 
$m_{1}^{\rho },m_{2}^{\rho },m_{3}^{\rho },...$

\begin{theorem}
(Equivariant Localization Theorem)\label{EquLocalizationThm}For each
irreducible representation $\rho $ of $G$ and for each fixed $N>0$, there
exists $c>0$ and $s_{0}>0$ such that for any $s>s_{0}$ and any $j\leq N$,
the interval $\left( \mu _{j}^{\rho }-cs^{-1/5},\mu _{j}^{\rho
}+cs^{-1/5}\right) $ contains exactly $m_{j}^{\rho }$ eigenvalues of $%
H_{s}^{\rho }$. Furthermore, all the eigenvalues of $H_{s}^{\rho }$
contained in $\left( -\infty ,\mu _{N}^{\rho }+cs^{-1/5}\right) $ are
contained in 
\begin{equation*}
\bigcup_{j\leq N}\left( \mu _{j}^{\rho }-cs^{-1/5},\mu _{j}^{\rho
}+cs^{-1/5}\right) .
\end{equation*}
\end{theorem}

\begin{proof}
We show how to generalize Theorem 1.1 in \cite{Shu1} to the equivariant
setting. We identify the parameter $s$ in our theorem with $\frac{1}{h}$ in 
\cite{Shu1}.

To obtain an upper bound for the eigenvalues of $H_{s}^{\rho }$ (or a lower
bound on the spectral counting function of $H_{s}^{\rho }$), we use
eigensections of the model operator $K^{\rho }$ to produce test sections for 
$H_{s}^{\rho }$. Suppose that $\psi $ is an eigensection of $K_{\overline{x}%
}^{\rho }:\Gamma \left( \mathbb{R}^{2n},E_{\overline{x}}\right) ^{\rho
}\rightarrow \Gamma \left( \mathbb{R}^{2n},E_{\overline{x}}\right) ^{\rho }$
corresponding to the eigenvalue $\lambda $. Let $J\in C_{0}^{\infty }\left( 
\mathbb{R}^{2n}\right) $ be a radial function defined such that $0\leq J\leq
1$ , $J\left( x\right) =1$ if $\left| x\right| \leq 1$, $J\left( x\right) =0$
if $\left| x\right| \geq 2$. For any $s>0$, let $J^{\left( s\right) }\left(
x\right) =J\left( s^{2/5}x\right) $. Then the section 
\begin{equation*}
\phi \left( x\right) =J^{\left( s\right) }\left( x\right) s^{n/2}\psi \left(
s^{1/2}x\right)
\end{equation*}
is in $\Gamma \left( \mathbb{R}^{2n},E_{\overline{x}}\right) ^{\rho }$ as
well, because $J^{\left( s\right) }$ is $G$-invariant. We produce a
corresponding element $\widetilde{\phi }\in \Gamma \left( M,E\right) ^{\rho
} $ that has support in a small neighborhood $U$ of $\overline{x}$ as
follows. Let $\gamma $ be the unit speed geodesic from $\overline{x}$ to $%
p\in U$, let $x_{p}$ be the geodesic normal coordinates of $p$, and let $%
P_{\gamma }:E_{\overline{x}}\rightarrow E_{p}$ denote parallel translation
along $\gamma $. We define 
\begin{equation*}
\widetilde{\phi }\left( p\right) =P_{\gamma }\phi \left( x_{p}\right) .
\end{equation*}
Clearly, $\widetilde{\phi }\in \Gamma \left( M,E\right) $. Because the
connection on $E$ is equivariant, parallel translation commutes with the
action of $G$, and $\widetilde{\phi }\in \Gamma \left( M,E\right) ^{\rho }$.
This specific trivialization of $E$ produces test sections that can be used
as in \cite{Shu1} to obtain the upper bounds for the eigenvalues of $%
H_{s}^{\rho }$. We denote $\Phi :\Gamma \left( \mathbb{R}^{2n},E_{\overline{x%
} }\right) ^{\rho }\rightarrow \Gamma \left( U,E\right) ^{\rho }$ to be the
trivialization $\phi \rightarrow \widetilde{\phi }$.

To obtain a lower bound on the eigenvalues of $H_{s}^{\rho }$ (or an upper
bound on the spectral counting function of $H_{s}^{\rho }$), we proceed
exactly as in \cite{Shu1}. The functions in the partition of unity are
chosen so that those corresponding to neighborhoods of singular points are
radial; then the partition of unity will consist of invariant functions.
Next, the IMS localization formula allows us to localize to these small
neighborhoods, comparing the operators $\Phi ^{-1}H_{s}^{\rho }\Phi $ and $%
K^{\rho }$.
\end{proof}

\section{ Analysis of Equivariant Perturbations\label{EquivariantWitten}}

\medskip In this section, we are going to apply Theorem \ref%
{EquLocalizationThm} to the following situation. Let $G=T^{m}\cong \mathbb{R}%
^{m}\diagup 2\pi \mathbb{Z}^{m}$ act on the right by isometries on a closed,
oriented Riemannian manifold $M$ of dimension $2n$. Let $D^{+}:\Gamma \left(
M,E^{+}\right) \rightarrow \Gamma \left( M,E^{-}\right) $ be a first-order, $%
G$-equivariant, transversally elliptic operator, where $E^{+}$ and $E^{-}$
are $G$-equivariant Hermitian vector bundles of rank $2r$ over $M$. Let $%
E=E^{+}\oplus E^{-}$, and let $D:\Gamma \left( M,E\right) \rightarrow \Gamma
\left( M,E\right) $ denote the operator $\left( D^{+},\left( D^{+}\right)
^{\ast }\right) $, where $\ast $ denotes the adjoint.

Consider the following family of operators, depending on a real parameter $%
s: $%
\begin{equation*}
D_{s}=D+sZ,
\end{equation*}%
where $Z$ has the following properties:

\begin{enumerate}
\item $Z:E^{\pm }\rightarrow E^{\mp }$ is a smooth, self-adjoint,
equivariant bundle map that is nonsingular away from a finite number of
points of $M$.

\item For each irreducible representation $\rho $, the restriction of $DZ+ZD$
to $\Gamma \left( M,E\right) ^{\rho }$ is a bundle map.

\item In a small neighborhood $U_{\overline{x}}$ of each singular point $%
\overline{x}$ of $Z$, we assume that $E^{\pm }$ has the structure of an
equivariant Clifford bundle (with equivariant Clifford connection $\nabla $)
and that $D$ is a(n equivariant) Dirac operator near these points (see \cite%
{B-G-V} ).

\item In $U_{\overline{x}}$, the operator $D+sZ$ has the following explicit
form. We require 
\begin{equation*}
Z=c\left( iV\right) ,\text{ so that }D_{s}=D+sc\left( iV\right) ,
\end{equation*}%
where $V$ is a vector field induced from some element $\mathbf{v}_{\overline{%
x}}$ of the Lie algebra $\mathfrak{g}$ of the torus $G$ such that the
closure of $\left\{ \exp \left( t\mathbf{v}_{\overline{x}}\right) \,|\,t\in 
\mathbb{R}\right\} $ is the entire torus $T^{m}$, and where $c\left(
iV\right) $ denotes Clifford multiplication by $iV$.
\end{enumerate}

For example, if $D$ is a Dirac operator on sections of a Clifford bundle and 
$V$ is a global Killing vector field with isolated fixed points that induces
an infinitesimal isometry of the bundle, then the operator $Z=c\left(
iV\right) $ satisfies the conditions above (see Lemma \ref{DZ+ZD Lemma}),
where the torus group is the closure of the flow of $V$ in the isometry
group of $M$. For a case of a transversally elliptic operator and
perturbation $Z$, see Example \ref{transvDiracSphere}.

The proof of the next lemma can be found in the Appendix.

\begin{lemma}
\label{DZ+ZD Lemma}In the notation above, for any vector field $V$, if $D$
is a Dirac operator, 
\begin{equation*}
\left( D_{s}\right) ^{2}=D^{2}+s\left( -2i\nabla _{V}-i\mathrm{div}\left(
V\right) +ic\left( d\left( V^{\ast }\right) \right) \right) +s^{2}\left\vert
V\right\vert ^{2}.
\end{equation*}%
Here $V^{\ast }$ is the one-form dual to the vector field $V$, and $c\left(
\alpha \wedge \beta \right) :=c\left( \alpha \right) c\left( \beta \right) $
for orthogonal covectors $\alpha $ and $\beta $.
\end{lemma}

In what follows, we need to define the Lie derivative of a section of $E$.
Since $G$ acts on $M$ on the right and since $E$ is $G$-equivariant, the
bundle $E$ is endowed with the lifted left action $F_{g}:E_{x}\rightarrow
E_{xg}$ on $E$ for each $g\in G$.

\begin{definition}
The induced action $\psi _{g}$ of $g\in G$ on the a section $u\in \Gamma
\left( M,E\right) $ is%
\begin{equation*}
\left( \psi _{g}u\right) \left( x\right) =F_{g^{-1}}\left( u\left( xg\right)
\right) ,
\end{equation*}%
and the action satisfies%
\begin{equation*}
\psi _{gh}=\psi _{h}\circ \psi _{g}
\end{equation*}%
for all $g,h\in G$.
\end{definition}

\begin{definition}
The Lie derivative $\mathcal{L}_{V}u$ of a section $u\in $ $\Gamma \left(
M,E\right) $ in direction $V$ (as above, the vector field induced from $%
\mathbf{v}\in \mathfrak{g}$) is%
\begin{equation*}
\left( \mathcal{L}_{V}u\right) \left( x\right) =\left. \frac{d}{dt}\left[
F_{\exp \left( -t\mathbf{v}\right) }\left( u\left( x\exp \left( t\mathbf{v}%
\right) \right) \right) \right] \right\vert _{t=0}.
\end{equation*}
\end{definition}

With this definition, $\mathcal{L}_{V}$ satisfies the usual properties of
Lie derivative on tensors. For example, the standard induced action of a Lie
group on vector fields and forms gives the ordinary Lie derivative. The
following lemma is standard.

\begin{lemma}
If $V$ is an infinitesimal isometry, then the operator $A_{V}=\nabla _{V}-%
\mathcal{L}_{V}$ is a skew-Hermitian endomorphism of $E$.
\end{lemma}

\begin{example}
\label{SpinorDZ+ZD}Let $V$ be a Killing field generating an action by
isometries on a Riemannian manifold $\left( M,g\right) .$ If $M$ is in
addition a spin manifold, then the action automatically lifts to the spinor
bundle $S$. If we let \emph{\ }$\mathcal{L}_{V}^{S}$ be the Lie derivative
of this action on the spinors, induced by the action on the frame bundle,
then 
\begin{equation*}
A_{V}=\nabla _{V}^{S}-\mathcal{L}_{V}^{S}=\frac{1}{4}c\left( d\left( V^{\ast
}\right) \right) .
\end{equation*}%
See the proof in the appendix.
\end{example}

\begin{corollary}
If $D$ is a Dirac operator, and if $D_{s}=D+sc\left( iV\right) $ with $V$ an
infinitesimal isometry as above, then 
\begin{equation}
\left( D_{s}\right) ^{2}=D^{2}+s\left( -2i\mathcal{L}_{V}-2iA_{V}+ic\left(
d\left( V^{\ast }\right) \right) \right) +s^{2}\left\vert V\right\vert ^{2}.
\label{DsSquaredEqu}
\end{equation}
\end{corollary}

\begin{proof}
Combine the two previous lemmas, and observe that for a Killing vector field 
$V$ we have $\mathrm{div}\left( V\right) =0.$
\end{proof}

\begin{remark}
\vspace{1pt}A similar computation was done by Bismut and explained in \cite[%
Chapter 8]{B-G-V} in the heat kernel proof of the Kirillov character
formula. In this computation, the endomorphism $-A_{X}$ is called the
\textquotedblleft moment\textquotedblright\ of the connection, used in the
context of frame bundles and bundles of forms.
\end{remark}

The nondegeneracy of $V$ at the zero $\overline{x}$ implies $\left\vert
V\left( x\right) \right\vert ^{2}\geq c\cdot d\left( x,\overline{x}\right)
^{2}$ for $x$ near $\overline{x}$\textbf{,} so the lemmas above imply that
the hypotheses of Theorem \ref{EquLocalizationThm} are satisfied for the
operator 
\begin{equation}
H_{s}=\frac{1}{s}\left( D_{s}\right) ^{2}.  \label{HsDefinition}
\end{equation}

\vspace{1pt}Fix $\rho :T^{m}\rightarrow \mathbb{C}$ to be a particular
irreducible unitary representation. Note that if we choose coordinates $%
\mathbf{\theta }=\left( \theta _{1},...,\theta _{m}\right) \in \left( 
\mathbb{R}\diagup 2\pi \mathbb{Z}\right) ^{m}$, then the representation has
the form 
\begin{equation}
\rho \left( \mathbf{\theta }\right) =e^{i\mathbf{b\cdot \theta }},
\label{rhoFormula}
\end{equation}%
where $\mathbf{b}=\left( b_{1},...,b_{m}\right) \in \mathbb{Z}^{m}$. Note
that the vector $\mathbf{b}$ depends on the choice of coordinates $\mathbf{%
\theta }$; for instance, if $\theta _{i}$ is replaced by $-\theta _{i}$,
then $b_{i}$ is replaced by $-b_{i}$. In what follows, the choice of
coordinates $\mathbf{\theta }$ will depend on $\mathbf{v}_{\overline{x}}$ .

Fix a singular point $\overline{x}$ of $Z$. We now describe the model
operator $K_{\overline{x}}^{\rho }$ and compute its eigenvalues.

We will use geodesic normal coordinates centered at $\overline{x}$. In these
coordinates, $\widetilde{A}$, the principal part of $A=D^{2}$ at $\overline{x%
}$, is the Euclidean Laplacian.

Now we compute 
\begin{equation*}
\widetilde{B}^{\rho }=\left. B\right\vert _{\Gamma \left( M,E\right) ^{\rho
}}\left( \overline{x}\right) =\left. \left( -2i\mathcal{L}%
_{V}-2iA_{V}+ic\left( d\left( V^{\ast }\right) \right) \right) \right\vert
_{\Gamma \left( M,E\right) ^{\rho }}\left( \overline{x}\right)
\end{equation*}%
The vector $\mathbf{v}_{\overline{x}}\in \mathfrak{g}$ generates a dense
flow $\mathbf{\theta }\left( t\right) $ on $T^{m}$ by the formula%
\begin{equation*}
\mathbf{\theta }\left( t\right) =\exp \left( t\mathbf{v}_{\overline{x}%
}\right) =t\mathbf{\tau }=\left( t\tau _{1},t\tau _{2},...,t\tau _{m}\right)
\in T^{m},
\end{equation*}%
where $\mathbf{\tau }=\left( \tau _{1},...,\tau _{m}\right) \in \mathbb{R}%
^{m}$. We choose the coordinates $\mathbf{\theta }$ so that the torus action
satisfies $\tau _{p}>0$ for $1\leq p\leq m$. Since the flow is dense, the
set $\left\{ \tau _{1},...,\tau _{m}\right\} $ is linearly independent over $%
\mathbb{Q}$.

The representation $\rho $ and choice of coordinates $\mathbf{\theta }$
uniquely determine the vector $\mathbf{b}\in \mathbb{Z}^{m}$ as in formula (%
\ref{rhoFormula}). If $u\in \Gamma \left( M,E\right) $ is of type $\rho $,
then near $\overline{x}$ we have 
\begin{equation*}
\mathcal{L}_{V}u=i\left( \mathbf{b}\cdot \mathbf{\tau }\right) u.
\end{equation*}

The action of $\mathbf{\theta }\in T^{m}$ on a small neighborhood of the
point $\overline{x}$ can be transferred to the tangent space $T_{\overline{x}%
}M$ via conjugation with the exponential map; the induced action on $T_{%
\overline{x}}M$ is an isometry.

Choose orthonormal coordinates $\left( x_{1},y_{1},...,x_{n},y_{n}\right) $ $%
=\left( z_{1},...,z_{n}\right) \ $on $T_{\overline{x}}M\cong \mathbb{C}^{n}$
so that $\mathbf{\theta }\in T^{m}$ acts by 
\begin{equation}
\left( z_{1},...,z_{n}\right) \mathbf{\theta }=\left( e^{i\mathbf{k}%
_{1}\cdot \mathbf{\theta }}z_{1},...,e^{i\mathbf{k}_{n}\cdot \mathbf{\theta }%
}z_{n}\right) ,  \label{kFormula}
\end{equation}%
where each $\mathbf{k}_{l}=\left( k_{l1},...,k_{lm}\right) \in \mathbb{Z}%
^{m} $. We assume in addition that for each $l$, 
\begin{equation}
\kappa _{l}:=\mathbf{k}_{l}\cdot \mathbf{\tau }>0;  \label{orientationChoice}
\end{equation}%
otherwise, replace $x_{l}$ with $y_{l}$ and vice versa. Note that the
resulting coordinates will not necessarily have the same orientation as the
induced orientation that comes from the manifold $M$.

Next, we choose an Hermitian coordinates $\left( w_{1},...,w_{r}\right) $ of 
$E_{\overline{x}}$ so that the action of $\mathbf{\theta }\in T^{m}$ on $E_{%
\overline{x}}$ is given by 
\begin{equation}
F_{\mathbf{\theta }}\left( w_{1},...,w_{r}\right) =\left( e^{i\mathbf{a}%
_{1}\cdot \mathbf{\theta }}w_{1},...,e^{i\mathbf{a}_{r}\cdot \mathbf{\theta }%
}w_{r}\right) ,  \label{aFormula}
\end{equation}%
with $\mathbf{a}_{j}=\left( a_{j1},...,a_{jm}\right) \in $ $\mathbb{Z}^{m}$.
Further, we choose the basis of $E_{\overline{x}}=E_{\overline{x}}^{+}\oplus
E_{\overline{x}}^{-}$ so that the grading operator is diagonal in this
basis. (Note that the grading commutes with the group action, so we may do
this.)

We compute that 
\begin{eqnarray*}
V &=&\sum_{l=1}^{n}\kappa _{l}\partial _{\phi _{l}},~~\left\vert
V\right\vert ^{2}=\sum_{l=1}^{n}\kappa _{l}^{2}\left\vert z_{l}\right\vert
^{2} \\
V^{\ast } &=&\sum_{l=1}^{n}\kappa _{l}\left\vert z_{l}\right\vert ^{2}d\phi
_{l},~~dV^{\ast }=2\sum_{l=1}^{n}\kappa _{l}d\mathrm{vol}_{l}, \\
A_{V} &=&\nabla _{V}-\mathcal{L}_{V}=i\sum_{j=1}^{r}\left( \mathbf{a}%
_{j}\cdot \mathbf{\tau }\right) P_{j},
\end{eqnarray*}%
where $\partial _{\phi _{l}}$ is the angular vector field $x_{l}\partial
_{y_{l}}-y_{l}\partial _{x_{l}}$ , $d\mathrm{vol}_{l}=dx_{l}\wedge dy_{l}$,
and $P_{j}=$ projection onto the $w_{j}$ plane (i.e. $j$th coordinate plane)
in $E_{\overline{x}}$.

Observe that the operators $ic\left( d\mathrm{vol}_{l}\right) $ mutually
commute, commute with the chirality operator and with the group action, and
square to $1$. Since the operators $\left[ ic\left( d\mathrm{vol}_{l}\right) %
\right] $ commute with the group action, they commute with each $P_{j}$. Let 
$\varepsilon _{jl}\in \left\{ -1,1\right\} $ be defined by 
\begin{equation}
\varepsilon _{jl}P_{j}=ic\left( d\mathrm{vol}_{l}\right) P_{j}
\label{epsilonFormula}
\end{equation}

Using these calculations, we obtain the second term of the model operator $%
K_{\overline{x}}^{\rho }$: 
\begin{eqnarray*}
\widetilde{B}^{\rho } &=&\left. \left( -2i\mathcal{L}_{V}-2iA_{V}+ic\left(
d\left( V^{\ast }\right) \right) \right) \right\vert _{\Gamma \left(
M,E\right) ^{\rho }}\left( \overline{x}\right) \\
&=&2\mathbf{b}\cdot \mathbf{\tau }+2\sum_{j=1}^{r}\left( \mathbf{a}_{j}\cdot 
\mathbf{\tau }\right) P_{j}+2\sum_{l=1}^{n}\kappa _{l}\left[ ic\left( d%
\mathrm{vol}_{l}\right) \right]
\end{eqnarray*}

Finally, we must compute 
\begin{equation*}
\widetilde{C}=\text{the quadratic part of }\left\vert V\right\vert ^{2}\text{
at }\overline{x}=\sum_{l=1}^{n}\kappa _{l}^{2}\left\vert z_{l}\right\vert
^{2}.
\end{equation*}

Thus, the model operator relevant to Theorem \ref{EquLocalizationThm} is 
\begin{eqnarray}
K_{\overline{x}}^{\rho } &=&\sum_{l=1}^{n}\left( -\partial
_{x_{l}}^{2}-\partial _{y_{l}}^{2}\right) +\left( 2\mathbf{b}\cdot \mathbf{%
\tau }+2\sum_{j=1}^{r}\left( \mathbf{a}_{j}\cdot \mathbf{\tau }\right)
P_{j}+2\sum_{l=1}^{n}\kappa _{l}\left[ ic\left( d\mathrm{vol}_{l}\right) %
\right] \right) +\sum_{l=1}^{n}\kappa _{l}^{2}\left\vert z_{l}\right\vert
^{2}  \notag \\
&=&\sum_{j=1}^{r}\left[ \sum_{l=1}^{n}\left( -\partial _{x_{l}}^{2}-\partial
_{y_{l}}^{2}+\sum_{l=1}^{n}\kappa _{l}^{2}\left( x_{l}^{2}+y_{l}^{2}\right)
\right) \right.  \notag \\
&&\left. +\left( 2\mathbf{b}\cdot \mathbf{\tau }+2\sum_{j=1}^{r}\left( 
\mathbf{a}_{j}\cdot \mathbf{\tau }\right) +2\sum_{l=1}^{n}\kappa
_{l}\varepsilon _{jl}\right) \right] P_{j}  \label{modelRestrictedTo}
\end{eqnarray}%
It is well-known that for the sum of oscillators 
\begin{equation*}
\sum_{l=1}^{n}\left( -\partial _{x_{l}}^{2}-\partial
_{y_{l}}^{2}+\sum_{l=1}^{n}\kappa _{l}^{2}\left( x_{l}^{2}+y_{l}^{2}\right)
\right) ,
\end{equation*}%
the eigenvalues are the numbers (determined by an arbitrary $\mathbf{m}%
=\left( m_{1},...,m_{n}\right) \in \mathbb{Z}^{n}$ and $\mathbf{d}=\left(
d_{1},...,d_{n}\right) \in \left( \mathbb{Z}_{\geq 0}\right) ^{n}$) 
\begin{equation*}
\lambda _{\mathbf{m},\mathbf{d}}=2\sum_{l=1}^{n}\kappa _{l}\left( \left\vert
m_{l}\right\vert +2d_{l}+1\right) ~,
\end{equation*}%
corresponding to the scalar eigenfunctions 
\begin{equation*}
\phi _{\mathbf{m},\mathbf{d}}=\prod_{l=1}^{n}e^{-\frac{1}{2}r_{l}^{2}\kappa
_{l}}r_{l}^{\left\vert m_{l}\right\vert }e^{im_{l}\theta _{l}}\cdot
L_{d_{l},\left\vert m_{l}\right\vert }\left( r_{l}^{2}\kappa _{l}\right) ,
\end{equation*}%
where $\left( r_{l},\theta _{l}\right) $ are polar coordinates in the $%
\left( x_{l},y_{l}\right) $-plane and $L_{d_{l},\left\vert m_{l}\right\vert
}\left( r\right) $ is a generalized Laguerre polynomial of degree $d_{l}\geq
0$. In particular, $L_{0,\left\vert m_{l}\right\vert }\left( r\right) =1$
for all $m_{l}\in \mathbb{Z}$. It is well known that the set $\left\{ \left.
\phi _{\mathbf{m},\mathbf{d}}\right\vert \,\mathbf{m}\in \mathbb{Z}^{n},%
\mathbf{d}\in \left( \mathbb{Z}_{\geq 0}\right) ^{n}\right\} $ is a
orthogonal basis of $L^{2}\left( \mathbb{C}^{n}\right) $.

Next, we compute the action of $\mathbf{\theta }\in T^{m}$ on each $\phi _{%
\mathbf{m},\mathbf{d}}$: 
\begin{eqnarray*}
\phi _{\mathbf{m},\mathbf{d}}\left( z_{1},...,z_{n}\right) &\mapsto &\phi _{%
\mathbf{m},\mathbf{d}}\left( e^{i\mathbf{k}_{1}\cdot \mathbf{\theta }%
}z_{1},...,e^{i\mathbf{k}_{n}\cdot \mathbf{\theta }}z_{n}\right) \\
&=&\exp \left( i\sum_{j=1}^{n}m_{j}\mathbf{k}_{j}\cdot \mathbf{\theta }%
\right) \phi _{\mathbf{m},\mathbf{d}}\left( z_{1},...,z_{n}\right) .
\end{eqnarray*}

For each $j\in \left\{ 1,...,r\right\} $, let $e_{j}$ be a basis vector that
spans $P_{j}\left( E_{\overline{x}}\right) $. Since we wish to consider
sections of type $\rho $, first observe that since $\mathbf{\theta }\in
T^{m} $ acts on a section $u\in \Gamma \left( \mathbb{R}^{2n},\mathbb{C}%
\right) $ by $\psi _{\mathbf{\theta }}\left( u\right) \left( z\right) =F_{-%
\mathbf{\theta }}u\left( z\mathbf{\theta }\right) $, 
\begin{eqnarray*}
\psi _{\mathbf{\theta }}\left( \phi _{\mathbf{m},\mathbf{d}}e_{j}\right)
\left( z\right) &=&F_{-\mathbf{\theta }}\left( \phi _{\mathbf{m},\mathbf{d}%
}\left( \left( z\right) \mathbf{\theta }\right) e_{j}\right) =F_{-\mathbf{%
\theta }}\left( \exp \left( i\sum_{h=1}^{n}m_{h}\mathbf{k}_{h}\cdot \mathbf{%
\theta }\right) \phi _{\mathbf{m},\mathbf{d}}\left( z\right) e_{j}\right) \\
&=&\exp \left( i\left[ -\mathbf{a}_{j}+\sum_{h=1}^{n}m_{h}\mathbf{k}_{h}%
\right] \cdot \mathbf{\theta }\right) \phi _{\mathbf{m},\mathbf{d}}\left(
z\right) e_{j}
\end{eqnarray*}%
In order that $\phi _{\mathbf{m},\mathbf{d}}P_{j}\in \Gamma \left( \mathbb{R}%
^{2n},\mathbb{C}\right) ^{\rho }$, the following equation must be satisfied:

\begin{equation*}
-\mathbf{a}_{j}+\sum_{h=1}^{n}m_{h}\mathbf{k}_{h}=\mathbf{b}~.
\end{equation*}%
We note that there are many choices of the integers $m_{h}$, in general an
infinite number, that satisfy the equations above for given $\mathbf{a}_{j}$%
, $\mathbf{k}_{h}$ , and $\mathbf{b}$. Taking the dot product with $\mathbf{%
\tau }$, we have 
\begin{eqnarray}
-\mathbf{a}_{j}\cdot \mathbf{\tau }+\sum_{h=1}^{n}m_{h}\left( \mathbf{k}%
_{h}\cdot \mathbf{\tau }\right) &=&\mathbf{b}\cdot \mathbf{\tau },\text{ or}
\notag \\
\sum_{h=1}^{n}m_{h}\kappa _{h} &=&\mathbf{a}_{j}\cdot \mathbf{\tau }+\mathbf{%
b}\cdot \mathbf{\tau ~}.  \label{CoefficientCondition}
\end{eqnarray}%
The possible $\mathbf{m\in }\mathbb{Z}^{n}$ satisfying (\ref%
{CoefficientCondition}) are integer points in an $\left( n-1\right) $%
-dimensional plane, since the right hand side is fixed. From (\ref%
{modelRestrictedTo}), the restriction of $K_{\overline{x}}^{\rho }$ to a
section $\phi _{\mathbf{m},\mathbf{d}}e_{j}$ of type $\rho $ with a specific
choice of the $\mathbf{m}$ gives the formula 
\begin{eqnarray*}
K_{\overline{x}}^{\rho }\phi _{\mathbf{m},\mathbf{d}}e_{j} &=&\left(
2\sum_{l=1}^{n}\kappa _{l}\left( \left\vert m_{l}\right\vert
+2d_{l}+1\right) +2\mathbf{b}\cdot \mathbf{\tau }+2\mathbf{a}_{j}\cdot 
\mathbf{\tau }+2\sum_{l=1}^{n}\kappa _{l}\varepsilon _{jl}\right) \phi _{%
\mathbf{m}}e_{j} \\
&=&\left( 2\sum_{l=1}^{n}\kappa _{l}\left( \left\vert m_{l}\right\vert
+m_{l}+2d_{l}+1+\varepsilon _{jl}\right) \right) \phi _{\mathbf{m}}e_{j}~.
\end{eqnarray*}

\medskip We have proved the following theorem:

\begin{theorem}
\label{ModelSpectrumTheorem}\medskip The spectrum of $K_{\overline{x}}^{\rho
}$ is the set of real numbers of the form 
\begin{equation*}
\lambda =2\sum_{l=1}^{n}\kappa _{l}\left( \left\vert m_{l}\right\vert
+m_{l}+2d_{l}+1+\varepsilon _{jl}\right) ,
\end{equation*}%
where the multiplicity of the eigenvalue $\lambda $ is the number of pairs $%
\left( \mathbf{m},\mathbf{d}\right) \in \mathbb{Z}^{n}\times \left( \mathbb{Z%
}_{\geq 0}\right) ^{n}\mathbb{\,}$such that there exists $j\in \left\{
1,...,r\right\} $ such that
\end{theorem}

\begin{equation*}
\sum_{h=1}^{n}m_{h}\mathbf{k}_{h}=\,\mathbf{a}_{j}+\mathbf{b}\text{ ~and }%
2\sum_{l=1}^{n}\kappa _{l}\left( \left\vert m_{l}\right\vert
+m_{l}+2d_{l}+1+\varepsilon _{jl}\right) =\lambda .
\end{equation*}

\begin{remark}
\medskip Note that the multiplicities of the eigenvalues above are finite,
since 
\begin{equation*}
\lambda =2\sum_{l=1}^{n}\kappa _{l}\left( \left\vert m_{l}\right\vert
+2d_{l}+1+\varepsilon _{jl}\right) +\mathbf{b}\cdot \mathbf{\tau }+\mathbf{a}%
_{j}\cdot \mathbf{\tau },
\end{equation*}%
and the quantities $\left\vert m_{l}\right\vert $ and $d_{l}$ must be
bounded.
\end{remark}

Since

\begin{equation*}
K^{\rho }=\bigoplus_{\text{singular point }\overline{x}}K_{\overline{x}%
}^{\rho },
\end{equation*}%
the spectrum $\sigma \left( K^{\rho }\right) $ satisfies 
\begin{equation*}
\sigma \left( K^{\rho }\right) =\bigcup_{\text{singular point }\overline{x}%
}\sigma \left( K_{\overline{x}}^{\rho }\right) .
\end{equation*}

\begin{remark}
\medskip Theorem \ref{EquLocalizationThm}, Theorem \ref{ModelSpectrumTheorem}
, and Equation (\ref{HsDefinition}) imply that as $s\rightarrow \infty $,
the eigenvalues of $\frac{1}{s}\left( D+sZ\right) ^{2}$ restricted to
sections of type $\rho $ approach the eigenvalues $\lambda $ of the $K^{\rho
}$ as described above.
\end{remark}

\section{Applications to Equivariant Index Theory\label{EquIndexThmSection}}

We now apply Theorem \ref{EquLocalizationThm} and Theorem \ref%
{ModelSpectrumTheorem} to compute the index $\mathrm{ind}_{T^{m}}^{\rho
}\left( D\right) $ of $D$ restricted to sections of type $\rho $. Since the
equivariant index does not depend on continuous perturbations, the index of $%
D$ restricted to sections of type $\rho $ is 
\begin{eqnarray*}
\mathrm{ind}_{T^{m}}^{\rho }\left( D\right) &=&\mathrm{ind}_{T^{m}}^{\rho
}\left( D_{s}\right) \\
&=&\dim \ker \left( \left. \left( D_{s}\right) ^{2}\right\vert _{\Gamma
\left( M,E^{+}\right) ^{\rho }}\right) -\dim \ker \left( \left. \left(
D_{s}\right) ^{2}\right\vert _{\Gamma \left( M,E^{-}\right) ^{\rho }}\right)
.
\end{eqnarray*}%
We now calculate these kernels independently using Theorem \ref%
{EquLocalizationThm} and Theorem \ref{ModelSpectrumTheorem}. The standard
argument implies the following lemma.

\begin{lemma}
The index satisfies 
\begin{equation*}
\mathrm{ind}_{T^{m}}^{\rho }\left( D\right) =\sum_{\overline{x}}\dim \ker K_{%
\overline{x}}^{\rho ,+}-\dim \ker K_{\overline{x}}^{\rho ,-},
\end{equation*}%
where the superscript $\pm $ refers to the restriction to $E_{\overline{x}%
}^{\pm }$.
\end{lemma}

Next, $\dim \ker \left( K_{\overline{x}}^{\rho }\right) $ is the number of
pairs $\left( \mathbf{m},\mathbf{d}\right) \in \mathbb{Z}^{n}\times \left( 
\mathbb{Z}_{\geq 0}\right) ^{n}\mathbb{\,}$such that there exists $j\in
\left\{ \text{ }1,...,r\right\} $ such that

\begin{equation*}
\sum_{h=1}^{n}m_{h}\mathbf{k}_{h}=\,\mathbf{a}_{j}+\mathbf{b}\text{ ~and }%
2\sum_{l=1}^{n}\kappa _{l}\left( \left\vert m_{l}\right\vert
+m_{l}+2d_{l}+1+\varepsilon _{jl}\right) =0\text{.}
\end{equation*}%
In this formula, the quantities $\mathbf{a}_{j}$, $\mathbf{k}_{h}$, $\mathbf{%
b}$, $\kappa _{h}$ all depend on the critical point $\overline{x}$. Since
each $\kappa _{h}$ is positive, $2\sum_{h=1}^{n}\kappa _{h}\left( \left\vert
m_{h}\right\vert +m_{h}+2d_{l}+1+\varepsilon _{jh}\right) =0$ if and only if
each $m_{h}$ is nonpositive, each $d_{l}$ is zero, and each $\varepsilon
_{jh}$ is $-1$ for $1\leq h\leq n$. Thus, we may express $\dim \ker \left(
K_{\overline{x}}^{\rho }\right) $ as 
\begin{multline*}
\#\Bigg\{\left. \mathbf{m}=\left( m_{1},...,m_{n}\right) \in \mathbb{Z}^{n}%
\mathbb{\,}\right\vert \text{ }m_{h}\leq 0\text{ and there exists }j\in
\left\{ \text{ }1,...,r\right\} \text{ such that} \\
\varepsilon _{jh}=-1\text{ for all }h,1\leq h\leq n~\text{and }%
\sum_{h=1}^{n}m_{h}\mathbf{k}_{h}=\,\mathbf{a}_{j}+\mathbf{b}\Bigg\}
\end{multline*}

The theorem below follows immediately. Recall that the $j^{\mathrm{th}}$
coordinate plane as in Formula (\ref{epsilonFormula}) is a subspace of $E_{%
\overline{x}}^{+}$ or of $E_{\overline{x}}^{-}$.

\begin{theorem}
(Transverse Index Theorem) \label{transverseGIndexTheorem}Let 
\begin{equation*}
k_{j}\left( \overline{x}\right) =\left\{ 
\begin{array}{ll}
\#\left\{ \left. \mathbf{m}\in \mathbb{Z}^{n}\mathbb{\,}\right\vert \text{ }%
m_{h}\leq 0~\text{and }\sum_{h=1}^{n}m_{h}\mathbf{k}_{h}=\,\mathbf{a}_{j}+%
\mathbf{b}\right\} ~ & \text{if }\varepsilon _{jh}=-1\text{ for all }h \\ 
0 & \text{otherwise}%
\end{array}%
\right.
\end{equation*}%
and let 
\begin{equation*}
\mathrm{sign}\left( j\right) =\pm 1,
\end{equation*}%
according to whether the $j^{\mathrm{th}}$ coordinate plane is in $E_{%
\overline{x}}^{\pm }$. Then 
\begin{equation*}
\mathrm{ind}_{T^{m}}^{\rho }\left( D\right) =\sum_{Z\left( \bar{x}\right)
=0}\sum_{j=1}^{r}\mathrm{sign}\left( j\right) k_{j}\left( \overline{x}%
\right) .
\end{equation*}
\end{theorem}

\section{Index multiplicities for the Lie group and its maximal torus. \label%
{torusindices}}

Suppose that $T^{m}$ is a maximal torus in a compact, connected Lie group $G$%
, and let $D$ be a $G$-equivariant, transversally elliptic, first order
differential operator that is also transversally elliptic with respect to
the $T^{m}$ action. Then there is a relationship between the multiplicities $%
\mathrm{ind}_{T^{m}}^{\rho }\left( D\right) $ and $\mathrm{ind}_{G}^{\mu
}\left( D\right) $, for a given irreducible representation $\mu $ of $G$. We
choose the coordinates $\mathbf{\theta }\in \mathbb{R}^{m}\diagup 2\pi 
\mathbb{Z}^{m}$ for the torus $T^{m}$. For any character $\xi _{\alpha }$ of
a (not necessarily irreducible) representation $\alpha $ of $G$, the
restriction of $\xi _{\alpha }$ of $G$ to $T^{m}$ yields a character of $%
T^{m}$. Let $\chi _{\mu }$ denote the character of a specific irreducible
unitary representation $\mu $, which has $L^{2}$ norm $1$ with respect to
the Haar measure. Since characters are class functions, the multiplicity $%
n_{\mu }^{\alpha }$ of $\mu $ in $\alpha $ is 
\begin{equation*}
n_{\mu }^{\alpha }=\int_{G}\xi _{\alpha }\left( g\right) \overline{\chi
_{\mu }\left( g\right) }\,dg=\int_{T^{m}}\widetilde{\xi _{\alpha }}\left( 
\mathbf{\theta }\right) \overline{\widetilde{\chi _{\mu }}\left( \mathbf{%
\theta }\right) \,}\,f\left( \mathbf{\theta }\right) \frac{d\mathbf{\theta }%
}{\left( 2\pi \right) ^{m}},
\end{equation*}%
where $f\left( \mathbf{\theta }\right) $ is the factor of the integrand from
the Weyl Integration Formula and $dg$ is the normalized, biinvariant volume
form on $G$, and where the tilde $\widetilde{\left( \cdot \right) }$ denotes
restriction to the torus. Since the characters $e^{i\mathbf{b}\cdot \mathbf{%
\theta }}$ of the torus form an orthonormal basis of $L^{2}\left(
T^{m}\right) $, we may write 
\begin{equation*}
\overline{\widetilde{\chi _{\mu }}\left( \mathbf{\theta }\right) \,}%
\,f\left( \mathbf{\theta }\right) =\sum_{\mathbf{b}}\beta _{\mu }^{\mathbf{b}%
}e^{-i\mathbf{b}\cdot \mathbf{\theta }},
\end{equation*}%
where $\beta _{\mu }^{\mathbf{b}}$ are complex numbers depending only on the
irreducible representation $\mu $. Then 
\begin{eqnarray*}
\int_{G}\xi _{\alpha }\left( g\right) \overline{\chi _{\mu }\left( g\right) }%
\,dg &=&\sum_{\mathbf{b}}\beta _{\mu }^{\mathbf{b}}\widetilde{n_{\mathbf{b}%
}^{\alpha }},\text{ where} \\
\widetilde{\xi _{\alpha }}\left( \mathbf{\theta }\right) &=&\sum_{b}%
\widetilde{n_{\mathbf{b}}^{\alpha }}e^{i\mathbf{b}\cdot \mathbf{\theta }},
\end{eqnarray*}%
so that $\widetilde{n_{\mathbf{b}}^{\alpha }}$ is the multiplicity of $e^{i%
\mathbf{b}\cdot \mathbf{\theta }}$ in the restriction of $\alpha $ to $T^{m}$%
. Thus, the multiplicities of $G$-irreducible representations in $\alpha $
are determined by the multiplicities of the $T^{m}$-irreducible
representations of the restriction of $\alpha $ to a maximal torus. Thus, we
have 
\begin{equation*}
\mathrm{ind}_{G}^{\mu }\left( D\right) =\sum_{\mathbf{b}}\beta _{\mu }^{%
\mathbf{b}}\mathrm{ind}_{T^{m}}^{\rho _{\mathbf{b}}}\left( D\right) ,
\end{equation*}%
where $\rho _{\mathbf{b}}\left( \mathbf{\theta }\right) $ is multiplication
by $e^{i\mathbf{b}\cdot \mathbf{\theta }}$. and so the index multiplicities
for the Lie group $G$ are determined in a universal way from the
multiplicities of the maximal torus. Note that the formula above is valid
even if $D$ is not transversally elliptic with respect to $T$, as long as
the multiplicities of the representations of type $\rho _{\mathbf{b}}$ in $%
\ker D$ and $\ker D^{\ast }$ with $\beta _{\mu }^{\mathbf{b}}\neq 0$ are
finite.

We comment that this procedure is a consequence of the following. If a
vector space is a unitary representation space of a compact, connected Lie
group $G$, it is also a representation space of the maximal torus $T$. This
vector space may be decomposed into irreducible representation spaces of $G$
or into irreducible representation spaces of $T$. If the multiplicities of
all of these irreducible representations are finite, then the $G$%
-multiplicities determine the $T$-multiplicities, and, surprisingly, the $T$%
-multiplicities determine the $G$-multiplicities.

\begin{example}
\label{SU(2)Example}Suppose that $G=SU\left( 2\right) $. We compute the
coefficients $\beta _{\mu }^{\mathbf{b}}$ for a given irreducible unitary
representation of $SU\left( 2\right) $. We follow \cite[pp. 84ff]{BrotD}.
Let $V_{n}$ be the space of homogeneous polynomials of degree $n$ in $%
z=\left( z_{1},z_{2}\right) \in \mathbb{C}^{2}$, and let%
\begin{equation*}
\mu _{n}:SU\left( 2\right) \rightarrow \mathrm{End}\left( V_{n}\right)
\end{equation*}%
be defined for $g\in SU\left( 2\right) $ by $\left( \mu _{n}\left( g\right)
P\right) \left( z\right) =P\left( zg\right) $. These are precisely the
irreducible unitary representations of $SU\left( 2\right) $. Let $T<SU\left(
2\right) $ be the maximal torus defined as 
\begin{equation*}
T=\left\{ \left. E\left( t\right) =\left( 
\begin{array}{cc}
e^{it} & 0 \\ 
0 & e^{-it}%
\end{array}%
\right) ~\right\vert ~t\in \mathbb{R}\right\} .
\end{equation*}%
The character $\chi _{n}$ of $\mu _{n}$ satisfies%
\begin{eqnarray*}
\chi _{n}\left( E\left( t\right) \right) &=&\widetilde{\chi _{\mu }}\left(
E\left( t\right) \right) =\sum_{k=0}^{n}e^{i\left( n-2k\right) t} \\
&=&\frac{\sin \left( \left( n+1\right) t\right) }{\sin \left( t\right) }%
\text{ for }t\notin \pi \mathbb{Z}.
\end{eqnarray*}%
For any class function $\omega $ on $SU\left( 2\right) $ we have the Weyl
integration formula%
\begin{equation*}
\int_{SU\left( 2\right) }\omega \left( g\right) ~dg=\int_{0}^{2\pi }\omega
\left( E\left( t\right) \right) ~\left( 2\sin ^{2}\left( t\right) \right) ~%
\frac{dt}{2\pi },
\end{equation*}%
so that the function $f$ from this section is defined by%
\begin{equation*}
f\left( t\right) =2\sin ^{2}\left( t\right)
\end{equation*}%
For generic $t$, 
\begin{eqnarray*}
\overline{\widetilde{\chi _{n}}\left( E\left( t\right) \right) \,}\,f\left(
t\right) &=&\sum_{b\in \mathbb{Z}}\beta _{n}^{b}e^{-ib\theta },\text{ or} \\
\frac{\sin \left( \left( n+1\right) t\right) }{\sin \left( t\right) }\left(
2\sin ^{2}\left( t\right) \right) &=&\sum_{b\in \mathbb{Z}}\beta
_{n}^{b}e^{-ib\theta }.
\end{eqnarray*}%
The left hand side is 
\begin{equation*}
2\left( \frac{e^{i\left( n+1\right) t}-e^{-i\left( n+1\right) t}}{2i}\right)
\left( \frac{e^{it}-e^{-it}}{2i}\right) =\allowbreak \frac{1}{2}e^{-int}-%
\frac{1}{2}e^{i\left( n+2\right) t}-\allowbreak \frac{1}{2}e^{-i\left(
n+2\right) t}+\frac{1}{2}e^{int}\allowbreak ,
\end{equation*}%
so that%
\begin{equation*}
\beta _{n}^{b}=\left\{ 
\begin{array}{ll}
\frac{1}{2} & \text{if }b=n\text{ or }-n \\ 
-\frac{1}{2} & \text{if }b=n+2\text{ or }-n-2 \\ 
0 & \text{otherwise.}%
\end{array}%
\right.
\end{equation*}%
Thus, if $D$ is an $SU\left( 2\right) $-equivariant, transversally elliptic,
first order differential operator on a closed manifold such that $D$ is also
transversally elliptic with respect to the circle action given by
restricting the $SU\left( 2\right) $ action to $T$, then 
\begin{eqnarray}
&&\mathrm{ind}_{SU\left( 2\right) }^{\mu _{n}}\left( D\right)  \notag \\
&=&\frac{1}{2}\left( \mathrm{ind}_{T}^{\rho _{n}}\left( D\right) +\mathrm{ind%
}_{T}^{\rho _{-n}}\left( D\right) -\mathrm{ind}_{T}^{\rho _{n+2}}\left(
D\right) -\mathrm{ind}_{T}^{\rho _{-n-2}}\left( D\right) \right) ,
\label{SU(2)mult formula}
\end{eqnarray}%
where the representation $\rho _{k}$ satisfies 
\begin{equation*}
\rho _{k}\left( E\left( t\right) \right) =\text{multiplication by }e^{ikt}.
\end{equation*}
\end{example}

\section{Examples\label{ExampleSection}}

\subsection{Signature and de Rham operators, torus action.}

Let $M$ be a Riemannian manifold of dimension $2n$ endowed with an isometric
action of $T^{m}$. Let $\overline{x}$ be an isolated fixed point of this
action. There exist geodesic normal coordinates $\left(
z_{1},...,z_{n}\right) \in \mathbb{C}^{n}$ and coordinates $\mathbf{\theta }%
\in \mathbb{R}^{m}\diagup 2\pi \mathbb{Z}^{m}$ for $T^{m}$ such that $%
\overline{x}=\left( 0,...,0\right) $ and the action of $\mathbf{\theta }$ is
expressed using the vectors $\mathbf{k}_{1}=\left( k_{11},...,k_{1m}\right)
,...,\mathbf{k}_{n}=\left( k_{n1},...,k_{nm}\right) \in \mathbb{Z}^{m}$ as
follows: 
\begin{equation*}
\left( z_{1},...,z_{n}\right) \mapsto \left( e^{i\mathbf{k}_{1}\cdot \mathbf{%
\theta }}z_{1},...,e^{i\mathbf{k}_{n}\cdot \mathbf{\theta }}z_{n}\right) .
\end{equation*}%
Let $\mathbf{v}$ be an element of the Lie algebra of $T^{m}$ such that%
\begin{equation*}
\mathbf{\theta }\left( t\right) =\exp \left( t\mathbf{v}\right) =\left(
t\tau _{1},t\tau _{2},...,t\tau _{m}\right) =t\mathbf{\tau }\in T^{m}
\end{equation*}%
generates a dense flow in $T^{m}$, so that the set $\left\{ \tau
_{1},...,\tau _{m}\right\} $ must be linearly independent over $\mathbb{Q}$,
and such that each $\tau _{p}>0$, as in Section \ref{EquivariantWitten}. Let 
$V$ be the vector field on $M$ generated by this action. Consider the
operator $d+d^{\ast }$ on forms $\Gamma \left( M,\Lambda ^{\ast }T^{\ast
}M\right) $, and let $c:T^{\ast }M\rightarrow \mathrm{End}\left( \Lambda
^{\ast }T^{\ast }M\right) $ denote the standard Clifford action by cotangent
vectors, so that the Dirac operator is 
\begin{equation*}
D=d+d^{\ast }=c\circ \nabla ,
\end{equation*}%
where $\nabla $ is the Levi-Civita connection on forms.

From (\ref{orientationChoice}), we have 
\begin{equation*}
\kappa _{q}=\mathbf{k}_{q}\cdot \mathbf{\tau }>0,\text{ }1\leq q\leq n;
\end{equation*}%
otherwise the orientation of the $q^{\mathrm{th}}$ plane needs to be
reversed. Assume that this has been done.

Next we compute the numbers $\varepsilon _{jq}$ from formula (\ref%
{epsilonFormula}). We consider the Hermitian operators 
\begin{equation*}
ic\left( d\mathrm{vol}_{q}\right) =ic\left( dx_{q}\right) c\left(
dy_{q}\right)
\end{equation*}%
evaluated at the fixed point. The vector space $E_{\overline{x}}=\Lambda
^{\ast }T_{\overline{x}}^{\ast }M$ consists of forms 
\begin{equation*}
\left( A_{1}+B_{1}dx_{1}+C_{1}dy_{1}+D_{1}dx_{1}\wedge dy_{1}\right) \wedge
...\wedge \left( A_{n}+B_{n}dx_{n}+C_{n}dy_{n}+D_{n}dx_{n}\wedge
dy_{n}\right) ,
\end{equation*}%
where $A_{q},B_{q},C_{q},D_{q}\in \mathbb{C}$ for each $q$. Observe that $%
ic\left( dx_{q}\right) c\left( dy_{q}\right) $ acts only on the $q^{\mathrm{%
th}}$ component of the wedge product above, and 
\begin{eqnarray*}
&&ic\left( dx_{q}\right) c\left( dy_{q}\right) \left(
A_{q}+B_{q}dx_{q}+C_{q}dy_{q}+D_{q}dx_{q}\wedge dy_{q}\right) \\
&=&i\left( A_{q}dx_{q}\wedge dy_{q}+B_{q}dy_{q}-C_{q}dx_{q}-D_{q}\right)
\end{eqnarray*}%
Hence, the $q^{\mathrm{th}}$ components of eigenspaces of $ic\left(
dx_{q}\right) c\left( dy_{q}\right) $ are 
\begin{eqnarray*}
E_{\pm 1} &=&\mathrm{span}\left\{ 1\pm idx_{q}\wedge dy_{q},dx_{q}\pm
idy_{q}\right\} \\
&=&\left\{ 
\begin{array}{ll}
\mathrm{span}\left\{ 1+\frac{1}{2}d\overline{z_{q}}\wedge
dz_{q},dz_{q}\right\} & \text{for eigenvalue }+1 \\ 
\mathrm{span}\left\{ 1-\frac{1}{2}d\overline{z_{q}}\wedge dz_{q},d\overline{%
z_{q}}\right\} & \text{for eigenvalue }-1%
\end{array}%
\right.
\end{eqnarray*}

Henceforth we choose a basis of the $4^{n}$-dimensional space $E_{\overline{x%
}}$ as follows. Let 
\begin{equation*}
\omega _{1q}=1+\frac{1}{2}d\overline{z_{q}}\wedge dz_{q};~\omega _{2q}=1-%
\frac{1}{2}d\overline{z_{q}}\wedge dz_{q};~\omega _{3q}=dz_{q};~\omega
_{4q}=d\overline{z_{q}}
\end{equation*}%
Then for each $\mathbf{i=}\left( i_{1},...,i_{n}\right) \in \left\{
1,2,3,4\right\} ^{n}$, we have the basis element $\omega _{\mathbf{i}}$
defined by 
\begin{equation*}
\omega _{\mathbf{i}}=\omega _{i_{1}1}\wedge ...\wedge \omega _{i_{n}n},
\end{equation*}%
and $\left\{ \left. \omega _{\mathbf{i}}\,\right\vert \,\mathbf{i}\in
\left\{ 1,2,3,4\right\} ^{n}\right\} $ forms a basis of $E_{\overline{x}}$.
Note that an element $\mathbf{\theta }\in T^{m}$ acts on 
\begin{equation*}
\left( \widetilde{A}_{1}\omega _{11}+\widetilde{B}_{1}\omega _{21}+%
\widetilde{C}_{1}\omega _{31}+\widetilde{D}_{1}\omega _{41}\right) \wedge
...\wedge \left( \widetilde{A}_{n}\omega _{1n}+\widetilde{B}_{n}\omega _{2n}+%
\widetilde{C}_{n}\omega _{3n}+\widetilde{D}_{n}\omega _{4n}\right)
\end{equation*}%
via 
\begin{equation*}
\widetilde{A}_{q}\omega _{1q}+\widetilde{B}_{q}\omega _{2q}+\widetilde{C}%
_{q}\omega _{3q}+\widetilde{D}_{q}\omega _{4q}\mapsto \widetilde{A}%
_{q}\omega _{1q}+\widetilde{B}_{q}\omega _{2q}+e^{-i\mathbf{k}_{q}\cdot 
\mathbf{\theta }}\ \widetilde{C}_{q}\omega _{3q}+e^{i\mathbf{k}_{q\;}\cdot 
\mathbf{\theta }}\widetilde{D}_{q}\omega _{4q}
\end{equation*}

The irreducible representation spaces are 1-dimensional spaces spanned by $%
\omega _{\mathbf{i}}$. The representation restricted to the \textrm{span}$%
\left\{ \omega _{\mathbf{i}}\right\} $ is $\rho \left( \mathbf{\theta }%
\right) =$ multiplication by 
\begin{equation*}
\prod_{q,i_{q}=3}e^{-i\mathbf{k}_{q}\cdot \mathbf{\theta }%
}\prod_{N,i_{N}=4}e^{+i\mathbf{k}_{N\;}\cdot \mathbf{\theta }},
\end{equation*}%
corresponding to the vector $\mathbf{a}_{\mathbf{i}}\in \mathbb{Z}^{m}$ from
(\ref{aFormula}) with 
\begin{equation*}
\mathbf{a}_{\mathbf{i}}=-\sum_{i_{q}=3}\mathbf{k}_{q}+\sum_{i_{q}=4}\mathbf{k%
}_{q}
\end{equation*}%
\newline

Next, we calculate the integers $\varepsilon _{\mathbf{i}q}$ from (\ref%
{epsilonFormula}):\newline
\begin{equation}
\varepsilon _{\mathbf{i}q}\omega _{\mathbf{i}}=ic\left( d\mathrm{vol}%
_{q}\right) \omega _{\mathbf{i}}=\left( -1\right) ^{i_{q}+1}\omega _{\mathbf{%
i}},\text{ so }\varepsilon _{\mathbf{i}q}=\left( -1\right) ^{i_{q}+1}.
\label{epsilonExampleFormula}
\end{equation}%
\newline
Theorem \ref{transverseGIndexTheorem} implies that the contribution to the
equivariant index of $D$ at $\overline{x}$ is an alternating sum of the
quantities%
\begin{equation*}
k_{\mathbf{i}}\left( \overline{x}\right) =\left\{ 
\begin{array}{ll}
\#\left\{ \left. \mathbf{m}\in \mathbb{Z}^{n}\mathbb{\,}\right\vert \text{ }%
m_{h}\leq 0~\text{and }\sum_{h=1}^{n}m_{h}\mathbf{k}_{h}=\,\mathbf{a}_{%
\mathbf{i}}+\mathbf{b}\right\} ~ & \text{if }\varepsilon _{\mathbf{i}h}=-1%
\text{ for all }h \\ 
0 & \text{otherwise}%
\end{array}%
\right.
\end{equation*}%
We only count those $k_{\mathbf{i}}\left( \overline{x}\right) $ with $i_{h}$
even for all $h\in \left\{ 1,...,n\right\} .$ For these, we have 
\begin{eqnarray}
k_{\mathbf{i}}\left( \overline{x}\right) &=&\#\left\{ \mathbf{m}\in \mathbb{Z%
}^{n}\,\left\vert \text{ }m_{h}\leq 0~\text{and }\mathbf{b}%
=\sum_{h=1}^{n}\left\{ 
\begin{array}{cc}
m_{h}\mathbf{k}_{h} & \text{if~}i_{h}=2 \\ 
\left( m_{h}-1\right) \mathbf{k}_{h} & \text{if~}i_{h}=4%
\end{array}%
\right. \right. \right\}  \notag \\
&=&\#\left\{ \mathbf{m}\in \mathbb{Z}^{n}\,\left\vert \text{ }m_{h}\leq 0~%
\text{and }\mathbf{b\cdot \tau }=\sum_{h=1}^{n}\left\{ 
\begin{array}{cc}
m_{h}\kappa _{h} & \text{if~}i_{h}=2 \\ 
\left( m_{h}-1\right) \kappa _{h} & \text{if~}i_{h}=4%
\end{array}%
\right. \right. \right\}  \label{linearcombinationb}
\end{eqnarray}%
The integer $\mathrm{sign}\left( \mathbf{i}\right) $ is $\pm 1$ according to
whether $\omega _{\mathbf{i}}\in E_{\overline{x}}^{+}$ or $E_{\overline{x}%
}^{-}$. For example, if $\mathbf{b}=\mathbf{0}$ (that is, we restrict to
invariant sections), then 
\begin{equation*}
k_{\mathbf{i}}\left( \overline{x}\right) =\left\{ 
\begin{array}{ll}
1~ & \text{if }i_{h}=2\text{ for every }h\in \left\{ 1,...,n\right\} . \\ 
0~ & \text{otherwise}%
\end{array}%
\right.
\end{equation*}%
Therefore, if $\rho _{0}$ is the trivial representation, 
\begin{equation*}
\mathrm{ind}_{T^{m}}^{\rho _{0}}\left( D\right) =\sum_{V\left( \overline{x}%
\right) =0}\pm 1,
\end{equation*}%
where the sign is determined by $\omega _{21}\wedge ...\wedge \omega
_{2n}\in E_{\overline{x}}^{+}$ or $E_{\overline{x}}^{-}$. Since $\omega
_{21}\wedge ...\wedge \omega _{2n}$ is an even form, for the de Rham
operator we have 
\begin{equation*}
\text{Euler}\left( M\right) ^{\rho _{0}}=\mathrm{ind}_{T^{m}}^{\rho
_{0}}\left( D\right) =\text{number of singular points of }V.
\end{equation*}%
If we consider the signature operator, the chirality of $\omega _{21}\wedge
...\wedge \omega _{2n}$ is $\left( -1\right) ^{n}\mathrm{sign}\left( V,%
\overline{x}\right) $, where $\mathrm{sign}\left( V,\overline{x}\right) $ is 
$\pm 1$ according to whether the orientation of $M$ agrees with the
orientation of our chosen coordinates $\left( z_{1},...,z_{n}\right) $ ---
that is, whether the orientation on the tangent space $T_{\overline{x}}M$
agrees with that induced from $V$. We write 
\begin{equation*}
\text{Signature}\left( M\right) ^{\rho _{0}}=\mathrm{ind}_{T^{m}}^{\rho
_{0}}\left( D\right) =\left( -1\right) ^{n}\sum_{V\left( \overline{x}\right)
=0}\mathrm{sign}\left( V,\overline{x}\right) .
\end{equation*}%
Note, that the kernel of $D$ consists of harmonic forms, which are always
invariant under isometric actions of connected Lie groups, so that in fact
the formulas above yield results about the Euler characteristic $\chi \left(
M\right) $ and signature: 
\begin{eqnarray}
\chi \left( M\right) &=&\text{number of singular points of }V
\label{Euler characteristic} \\
\text{Signature}\left( M\right) &=&\left( -1\right) ^{n}\sum_{V\left( 
\overline{x}\right) =0}\mathrm{sign}\left( V,\overline{x}\right) .
\label{Signature}
\end{eqnarray}

If our representation is not trivial, then $b_{p}\neq 0$ for some $p$, and 
\begin{equation*}
\mathbf{b}\cdot \mathbf{\tau }=-\sum_{h=1}^{n}c_{h}\kappa _{h}
\end{equation*}%
with integers $c_{h}\geq 0$ for all $h$ (see (\ref{linearcombinationb})). If 
$c_{h}>0$, both $i_{h}=2$ and $i_{h}=4$ yield positive values of $k_{\mathbf{%
i}}\left( \overline{x}\right) $. Let $A$ be a subset of $\left\{
1,...,n\right\} $, and let 
\begin{gather}
N\left( A,\mathbf{b},\overline{x}\right) =\text{\# of ways to write }\mathbf{%
b}\cdot \mathbf{\tau }=-\sum_{h\in A}c_{h}\kappa _{h}  \notag \\
\text{with }c_{h}\in \mathbb{Z}_{>0}\text{ for all }h\in A.
\label{chFormula}
\end{gather}%
Further, let 
\begin{equation*}
S\left( A,\overline{x}\right) =\sum_{\mathbf{i}\in \mathbf{I}_{A}}\mathrm{\
sign}\left( \mathbf{i}\right) ,
\end{equation*}%
where the sum is taken over the set $\mathbf{I}_{A}$ of all multi-indices $%
\mathbf{i}=\left( i_{1},...,i_{n}\right) $ such that 
\begin{equation*}
i_{h}=\left\{ 
\begin{array}{ll}
2\text{ or }4~ & \text{if }h\in A. \\ 
2~ & \text{otherwise}%
\end{array}%
\right. ,
\end{equation*}%
and where 
\begin{equation*}
\mathrm{sign}\left( \mathbf{i}\right) =\left\{ 
\begin{array}{ll}
1 & \text{if }\omega _{\mathbf{i}}\in E_{\overline{x}}^{+} \\ 
-1~ & \text{if }\omega _{\mathbf{i}}\in E_{\overline{x}}^{-}%
\end{array}%
\right. ,
\end{equation*}%
Then the contribution of the critical point $\overline{x}$ is%
\begin{equation*}
\sum_{\mathbf{i}}\mathrm{sign}\left( \mathbf{i}\right) k_{\mathbf{i}}\left( 
\overline{x}\right) =\sum_{A\subset \left\{ 1,...,n\right\} }N\left( A,%
\mathbf{b},\overline{x}\right) S\left( A,\overline{x}\right) ,
\end{equation*}%
and thus 
\begin{equation*}
\mathrm{ind}_{T^{m}}^{\rho _{\mathbf{b}}}\left( D\right) =\sum_{V\left( 
\overline{x}\right) =0}\sum_{A\subset \left\{ 1,...,n\right\} }N\left( A,%
\mathbf{b},\overline{x}\right) S\left( A,\overline{x}\right) .
\end{equation*}

Note that if $D$ is the de Rham operator and if $b_{p}\neq 0$ for some $p\in
\left\{ 1,...,m\right\} $, then $S\left( A,\overline{x}\right) =0$ for every
singular point $\overline{x}$ and every subset $A\subset \left\{
1,...,n\right\} $. The reason is that replacing $\omega _{2q}$ with $\omega
_{4q}$ or vice versa changes the parity of the form. Thus, 
\begin{equation}
\text{Euler}\left( M\right) ^{\rho _{\mathbf{b}}}=\mathrm{ind}_{T^{m}}^{\rho
_{\mathbf{b}}}\left( D\right) =0,  \label{b-Euler}
\end{equation}%
which agrees with the fact that the kernel of $D$ consists of invariant
forms.

If $D$ is the signature operator, then $\mathrm{sign}\left( \mathbf{i}%
\right) $ is the same for each of the indices such that $i_{h}=2$ or $4$ and
is $\mathrm{sign}\left( V,\overline{x}\right) $. Thus, if $b_{p}\neq 0$ for
some $p\in \left\{ 1,...,m\right\} $, 
\begin{eqnarray}
\text{Signature}\left( M\right) ^{\rho _{\mathbf{b}}} &=&\mathrm{ind}%
_{T^{m}}^{\rho _{\mathbf{b}}}\left( D\right)  \notag \\
&=&\sum_{V\left( \overline{x}\right) =0}\mathrm{sign}\left( V,\overline{x}%
\right) \sum_{A\subset \left\{ 1,...,n\right\} }2^{\left\vert A\right\vert
}N\left( A,\mathbf{b},\overline{x}\right) ,  \label{b-Signature}
\end{eqnarray}%
which is (surprisingly) zero because the harmonic forms are invariant.

\begin{proposition}
\label{killingProp}Given a Killing field $V$ on an $2n$-dimensional manifold
with isolated singularities,%
\begin{eqnarray*}
\chi \left( M\right) &=&\text{number of singular points of }V \\
\text{Signature}\left( M\right) &=&\left( -1\right) ^{n}\sum_{V\left( 
\overline{x}\right) =0}\mathrm{sign}\left( V,\overline{x}\right)
\end{eqnarray*}%
For any $\mathbf{b}\in \mathbb{Z}^{m}$, where $V$ generates a $T^{m}$
action, 
\begin{equation*}
\sum_{V\left( \overline{x}\right) =0}\mathrm{sign}\left( V,\overline{x}%
\right) \sum_{A\subset \left\{ 1,...,n\right\} }2^{\left\vert A\right\vert
}N\left( A,\mathbf{b},\overline{x}\right) =0.
\end{equation*}
\end{proposition}

\begin{remark}
The first two identities above are well known. The first identity is a
special case of the Hopf index theorem where the Hopf index of the vector
field is one at each singular point. The factor $\left( -1\right) ^{n}$ may
be removed from the second identity, since the signature is zero if $n$ is
odd. Both of the first two identities are particular cases of the
Atiyah-Bott fixed point formula. The third identity seems to be new.
\end{remark}

\begin{example}
A particular example of the calculations above is the following action of $%
T^{n+1}$ on complex projective space $\mathbb{C}P^{n}$. Consider homogeneous
coordinates $\left[ z_{1},...,z_{n+1}\right] $ with the standard metric, and
consider the family of isometries $\left[ z_{1},...,z_{n+1}\right] \mapsto %
\left[ e^{i\theta _{1}}z_{1},...,e^{i\theta _{n+1}}z_{n+1}\right] $, where $%
\theta =\left( \theta _{1},...,\theta _{n+1}\right) \in T^{n+1}$. Let 
\begin{equation*}
\theta \left( t\right) =\exp \left( t\mathbf{v}\right) =\left( t\tau
_{1},t\tau _{2},...,t\tau _{n+1}\right) \in T^{n+1},
\end{equation*}%
generates a dense flow in $T^{n+1}$ so that the set $\left\{ \tau
_{1},...,\tau _{n+1}\right\} $ is linearly independent over $\mathbb{Z}$ and 
$0<\tau _{1}<...<\tau _{n+1}$. Let $V$ be the vector field generated by this
action.

There are $n+1$ fixed points: $\left[ e_{1}\right] =\left[ 1,0,...,0\right] $%
, $\left[ e_{2}\right] =\left[ 0,1,0,...,0\right] $, ... , and $\left[
e_{n+1}\right] =\left[ 0,...,0,1\right] $. The homogeneous coordinates (in a
coordinate chart diffeomorphic to $\mathbb{C}^{n}$) near $\left[ e_{l}\right]
$ are $\left[ z_{1},...,z_{l-1},1,z_{l+1},...,z_{n+1}\right] $, and the
action in these coordinates is 
\begin{equation*}
\left[ z_{1},...,z_{l-1},1,z_{l+1},...,z_{n+1}\right] \mapsto \left[
e^{i\left( \theta _{1}-\theta _{l}\right) }z_{1},...,1,...,e^{i\left( \theta
_{n+1}-\theta _{l}\right) }z_{n+1}\right] .
\end{equation*}%
The numbers $k_{hp}$ in Formula \ref{kFormula} are 
\begin{equation*}
k_{hp}=\delta _{hp}-\delta _{pl},
\end{equation*}%
where $\delta _{hp}$ is the Kronecker delta. From (\ref{orientationChoice}), 
\begin{equation*}
\kappa _{q}=\mathbf{k}_{q}\cdot \mathbf{\tau }=\tau _{q}-\tau _{l}>0,~1\leq
q\leq n+1~\text{and~}q\neq l;
\end{equation*}%
Thus the orientation of the $q^{\mathrm{th}}$ plane needs to be reversed if $%
q<l$. Thus we have 
\begin{equation*}
\mathrm{sign}\left( V,\left[ e_{l}\right] \right) =\left( -1\right) ^{l+1}.
\end{equation*}

From Equations (\ref{Euler characteristic}), (\ref{Signature}), and (\ref%
{b-Euler}), we have 
\begin{eqnarray*}
\chi \left( \mathbb{C}P^{n}\right) &=&n+1 \\
\text{Signature}\left( \mathbb{C}P^{n}\right) &=&\sum_{l=1}^{n+1}\left(
-1\right) ^{l+1}=\left\{ 
\begin{array}{ll}
1 & \text{if }n\text{ is even} \\ 
0 & \text{if }n\text{ is odd}%
\end{array}%
\right. \\
\chi \left( \mathbb{C}P^{n}\right) ^{\rho _{\mathbf{b}}} &=&0
\end{eqnarray*}%
Formula (\ref{chFormula}) gives 
\begin{eqnarray*}
\mathbf{b}\cdot \mathbf{\tau } &=&\sum_{h=1}^{l-1}c_{h}\left( \tau _{h}-\tau
_{l}\right) -\sum_{h=l+1}^{n+1}c_{h}\left( \tau _{h}-\tau _{l}\right) ,\text{
so} \\
c_{1}
&=&b_{1};\,...\,;c_{l-1}=b_{l-1};c_{l+1}=-b_{l+1};\,...\,;c_{n+1}=-b_{n+1} \\
b_{l} &=&\sum_{h=l+1}^{n+1}c_{h}-\sum_{h=1}^{l-1}c_{h}.
\end{eqnarray*}%
Thus, the only irreducible representations $\rho _{\mathbf{b}}$ that could
have nonzero contributions from the singular point $\left[ e_{l}\right] $
are those which satisfy 
\begin{eqnarray*}
b_{1},...,b_{l-1} &\geq &0;\,b_{l+1},...,b_{n+1}\leq 0;\text{ and} \\
b_{l} &=&-\sum_{h\neq l}b_{h}.
\end{eqnarray*}%
Since the integers $c_{h}$ determine $\mathbf{b}$, we have 
\begin{equation*}
N\left( A,\mathbf{b},\left[ e_{l}\right] \right) =1.
\end{equation*}%
We give a specific example of a representation and the resulting formula; we
leave it to the reader to obtain a general formula that works for all
possible $\mathbf{b}$. Suppose 
\begin{equation*}
\left( b_{1},...,b_{13}\right) =\left( 0,1,0,0,76,0,0,0,0,0,-51,-24,-2\right)
\end{equation*}%
for the action of $T^{13}$ on $\mathbb{C}P^{12}$. Then, by the computations
above and formula (\ref{b-Signature}), we have 
\begin{eqnarray*}
\text{Signature}\left( M\right) ^{\rho _{\mathbf{b}}} &=&\sum_{V\left( 
\overline{x}\right) =0}\mathrm{sign}\left( V,\overline{x}\right)
\sum_{A\subset \left\{ 1,...,n\right\} }2^{\left\vert A\right\vert }N\left(
A,\mathbf{b},\overline{x}\right) \\
&=&\left( -1\right) ^{6}2^{4}+\left( -1\right) ^{7}2^{5}+\left( -1\right)
^{8}2^{5}+\left( -1\right) ^{9}2^{5}+ \\
&&+\left( -1\right) ^{10}2^{5}+\left( -1\right) ^{11}2^{5}+\left( -1\right)
^{12}2^{4} \\
&=&0.
\end{eqnarray*}
\end{example}

\subsection{An example of an $SU\left( 2\right) $-action}

As in Example \ref{SU(2)Example}, let $T<SU\left( 2\right) $ be the maximal
torus defined as 
\begin{equation*}
T=\left\{ \left. \left( 
\begin{array}{cc}
e^{it} & 0 \\ 
0 & e^{-it}%
\end{array}%
\right) ~\right\vert ~t\in \mathbb{R}\right\} .
\end{equation*}%
Consider the manifold $M=SU\left( 2\right) \diagup T$. We identify each $%
\left( \alpha ,\beta \right) \in S^{3}\subset \mathbb{C}$ with the
corresponding matrix $\left( 
\begin{array}{cc}
\alpha & -\overline{\beta } \\ 
\beta & \overline{\alpha }%
\end{array}%
\right) \in SU\left( 2\right) $. Each element of $M$ is an equivalence class
depending on $\left( \alpha ,\beta \right) \in S^{3}\subset \mathbb{C}$ : 
\begin{eqnarray*}
\left[ \left( \alpha ,\beta \right) \right] &=&\left\{ \left. \left( 
\begin{array}{cc}
\alpha & -\overline{\beta } \\ 
\beta & \overline{\alpha }%
\end{array}%
\right) \left( 
\begin{array}{cc}
e^{it} & 0 \\ 
0 & e^{-it}%
\end{array}%
\right) =\left( 
\begin{array}{cc}
\alpha e^{it} & -\overline{\beta }e^{-it} \\ 
\beta e^{it} & \overline{\alpha }e^{-it}%
\end{array}%
\right) ~\right\vert ~t\in \mathbb{R}\right\} \\
&=&\left[ \left( e^{it}\alpha ,e^{it}\beta \right) \right] .
\end{eqnarray*}%
We endow $M$ with the standard metric and the left $SU\left( 2\right) $%
-action, which is%
\begin{eqnarray*}
\left( z,w\right) \left[ \left( \alpha ,\beta \right) \right] &=&\left[
\left( 
\begin{array}{cc}
z & -\overline{w} \\ 
w & \overline{z}%
\end{array}%
\right) \left( 
\begin{array}{cc}
\alpha & -\overline{\beta } \\ 
\beta & \overline{\alpha }%
\end{array}%
\right) \right] \\
&=&\left[ \left( z\alpha -\overline{w}\beta ,w\alpha +\overline{z}\beta
\right) \right] .
\end{eqnarray*}%
Then $M=S^{2}$ (because $\left( \alpha ,\beta \right) \rightarrow \left[
\left( \alpha ,\beta \right) \right] $ is the Hopf fibration), so for
example its Euler characteristic is $2$ and its signature is zero. Consider
the de Rham operator $D=d+d^{\ast }$ on forms; this operator commutes with
the $SU\left( 2\right) $-action and with the even-odd grading. The kernel of 
$D$ is the set of harmonic forms, which consist of constants and constants
times the volume form. Both of these are invariant forms, so we have that
(in the notation of Example \ref{SU(2)Example}) 
\begin{equation*}
\mathrm{ind}_{SU\left( 2\right) }^{\mu _{n}}\left( D\right) =\left\{ 
\begin{array}{ll}
2~ & \text{if }n=0 \\ 
0 & \text{otherwise}%
\end{array}%
\right.
\end{equation*}%
Since the de Rham operator is transversally elliptic with respect to the
torus $T$ action, we also have%
\begin{equation*}
\mathrm{ind}_{SU\left( 2\right) }^{\mu _{n}}\left( D\right) =\frac{1}{2}%
\left( \mathrm{ind}_{T}^{\rho _{n}}\left( D\right) +\mathrm{ind}_{T}^{\rho
_{-n}}\left( D\right) -\mathrm{ind}_{T}^{\rho _{n+2}}\left( D\right) -%
\mathrm{ind}_{T}^{\rho _{-n-2}}\left( D\right) \right)
\end{equation*}%
by Equation \ref{SU(2)mult formula}. The $T$-action on $M$ is given by%
\begin{equation*}
\left( e^{is},0\right) \left[ \left( \alpha ,\beta \right) \right] =\left[
\left( e^{is}\alpha ,e^{-is}\beta \right) \right] =\left[ \left(
e^{2is}\alpha ,\beta \right) \right] .
\end{equation*}%
It has two fixed points, $\left[ \left( 1,0\right) \right] $ and $\left[
\left( 0,1\right) \right] $.

On one hand we know that the kernel of $D$ consists of $T$-invariant forms,
so that%
\begin{equation*}
\mathrm{ind}_{T}^{\rho _{n}}\left( D\right) =\left\{ 
\begin{array}{ll}
2~ & \text{if }n=0 \\ 
0 & \text{otherwise}%
\end{array}%
\right.
\end{equation*}%
for all $n\in \mathbb{Z}$. Thus, we have that if $n\geq 0$,%
\begin{eqnarray*}
\frac{1}{2}\left( \mathrm{ind}_{T}^{\rho _{n}}\left( D\right) +\mathrm{ind}%
_{T}^{\rho _{-n}}\left( D\right) -\mathrm{ind}_{T}^{\rho _{n+2}}\left(
D\right) -\mathrm{ind}_{T}^{\rho _{-n-2}}\left( D\right) \right) &=&\left\{ 
\begin{array}{ll}
\frac{1}{2}\left( 1+1-0-0\right) ~ & \text{if }n=0 \\ 
0 & \text{otherwise}%
\end{array}%
\right. \\
&=&\mathrm{ind}_{SU\left( 2\right) }^{\mu _{n}}\left( D\right) ,
\end{eqnarray*}%
as expected.

To compute the index of $D$ using Theorem \ref{transverseGIndexTheorem}, we
let the vector field $V$ be the infinitesimal generator of the action $\left[
\left( \alpha ,\beta \right) \right] \mapsto \left[ \left( e^{is}\alpha
,e^{-is}\beta \right) \right] $. Using the calculation in the previous
section, the index $\mathrm{ind}_{T}^{\rho _{0}}\left( D\right) $ is the
number of singular points, which is two, and all other indices are zero, as
expected.

Note that if we compute $\mathrm{ind}_{SU\left( 2\right) }^{\mu _{n}}$ or $%
\mathrm{ind}_{T}^{\rho _{n}}$ of the spin Dirac operator, we obtain zero for
all indices, by the Atiyah-Hirzebruch vanishing theorem.

Formula \ref{SU(2)mult formula} does not apply when the operator is
transversally elliptic with respect to the $SU\left( 2\right) $ action but
not to the $T$-action. For example, the zero operator%
\begin{equation*}
\mathbf{0}:\Gamma \left( SU\left( 2\right) \diagup T,\mathbb{C}\right)
\rightarrow \Gamma \left( SU\left( 2\right) \diagup T,\left\{ 0\right\}
\right) 
\end{equation*}%
is equivariant and transversally elliptic with respect to the $SU\left(
2\right) $ action and is equivariant but not transversally elliptic with
respect to the $T$ action. One may check that the $\mu _{n}$ part of $\ker 
\mathbf{0}$ is zero if $n$ is odd otherwise is the eigenspace of the
Laplacian on $S^{2}=SU\left( 2\right) \diagup T$ with eigenvalue $\frac{n}{2}
$, which occurs with multiplicity $1$. Thus,%
\begin{equation*}
\mathrm{ind}_{SU\left( 2\right) }^{\mu _{n}}\left( \mathbf{0}\right)
=\left\{ 
\begin{array}{ll}
0~ & \text{if }n\text{ is odd} \\ 
1 & \text{if }n\text{ is even}%
\end{array}%
\right. .
\end{equation*}%
for all $n\geq 0$. Note that in every irreducible representation space $\mu
_{n}$, the representations $\rho _{n},\rho _{n-2},...,\rho _{-n}$ of $T$
occur, each with multiplicity $1$. Thus the $\rho _{n}$ part of $\ker 
\mathbf{0}$ is%
\begin{equation*}
\mathrm{ind}_{SU\left( 2\right) }^{\rho _{n}}\left( \mathbf{0}\right)
=\left\{ 
\begin{array}{ll}
0~ & \text{if }n\text{ is odd} \\ 
\infty  & \text{if }n\text{ is even}%
\end{array}%
\right. .
\end{equation*}%
This demonstrates that Equation \ref{SU(2)mult formula} is valid only if the
corresponding indices $\mathrm{ind}_{T}^{\rho _{n}}\left( D\right) $ are
finite, which happens always when $D$ is $T$-transversally elliptic. Note
that the formula remains valid even if $D$ is not transversally elliptic if
the corresponding $\rho _{n}$ parts of the subspaces are finite dimensional,
as above in the case where $n$ is odd.

\subsection{A transversally elliptic operator on the sphere\label%
{transvDiracSphere}}

\subsubsection{The operator $D$ and its equivariant index}

Let $\alpha \in S^{1}$ act on 
\begin{equation*}
S^{2}=\left\{ \left( x,y,z\right) \in \mathbb{R}^{3}:x^{2}+y^{2}+z^{2}=1%
\right\} 
\end{equation*}
by a rotation of $2\alpha $ around the $z$-axis. Let $E$ be the trivial $%
\mathbb{C}^{2}$ bundle over $S^{2}$. Let $\alpha \in S^{1}$ act on $\left( 
\begin{array}{c}
w_{1} \\ 
w_{2}%
\end{array}%
\right) \in E$ by $F_{\alpha }\left( 
\begin{array}{c}
w_{1} \\ 
w_{2}%
\end{array}%
\right) =\left( 
\begin{array}{c}
e^{-i\alpha }w_{1} \\ 
e^{i\alpha }w_{2}%
\end{array}%
\right) $. We use the grading $E^{+}=\left\{ \left( 
\begin{array}{c}
f \\ 
0%
\end{array}%
\right) \right\} $, $E^{-}=\left\{ \left( 
\begin{array}{l}
0 \\ 
g%
\end{array}%
\right) \right\} $. Consider the transversally elliptic operator $D$ on
sections of $E$ over $S^{2}$. We will write it in two ways, using
rectangular coordinates $\left( x,y,z\right) $ or spherical coordinates $%
\left( \theta ,\phi \right) $ with $\phi $ the angle between the position
vector and the $z$-axis and $\theta $ the polar angle in the $xy$-plane. Let 
\textrm{Proj}$:T\mathbb{R}^{3}\rightarrow TS^{2}$ be the orthogonal
projection.

\begin{eqnarray*}
D &=&\left( 
\begin{array}{ll}
0 & -1 \\ 
1 & 0%
\end{array}%
\right) \text{\textrm{Proj}}\left( -z\frac{\partial }{\partial x}+x\frac{%
\partial }{\partial z}\right) +\left( 
\begin{array}{ll}
0 & i \\ 
i & 0%
\end{array}%
\right) \text{\textrm{Proj}}\left( -z\frac{\partial }{\partial y}+y\frac{%
\partial }{\partial z}\right) \\
&=&\left( 
\begin{array}{ll}
0 & e^{-i\theta } \\ 
-e^{i\theta } & 0%
\end{array}%
\right) \frac{\partial }{\partial \phi }+\cot \phi \left( 
\begin{array}{ll}
0 & -ie^{-i\theta } \\ 
-ie^{i\theta } & 0%
\end{array}%
\right) \frac{\partial }{\partial \theta }.
\end{eqnarray*}%
Note that this operator fails to be elliptic precisely at the equator $z=0$
(or $\phi =\frac{\pi }{2}$). It is an easy exercise to check that this
operator is $S^{1}$-equivariant and symmetric for the standard metric on $%
S^{2}$. It can be shown that these properties imply that $D$ is essentially
self-adjoint.

If $u=\left( 
\begin{array}{l}
u_{1} \\ 
u_{2}%
\end{array}%
\right) \in \ker D$, then in the upper hemisphere $z=\sqrt{1-x^{2}-y^{2}}$, 
\begin{eqnarray*}
D\left( 
\begin{array}{c}
u_{1} \\ 
0%
\end{array}%
\right) \left( x,y\right) &=&\left( 
\begin{array}{c}
0 \\ 
\left( -z\frac{\partial u_{1}}{\partial x}+x\frac{\partial u_{1}}{\partial z}%
\right) +i\left( -z\frac{\partial u_{1}}{\partial y}+y\frac{\partial u_{1}}{%
\partial z}\right)%
\end{array}%
\right) \\
&=&\left( 
\begin{array}{c}
0 \\ 
\left( -z\left( \frac{\partial }{\partial x}+i\frac{\partial }{\partial y}%
\right) \right) u_{1}\left( x,y\right)%
\end{array}%
\right) =\mathbf{0}.
\end{eqnarray*}%
Thus, $u_{1}$ must be holomorphic as a function of the coordinates $\left(
x,y\right) $. Similarly, $u_{2}$ must be antiholomorphic. The same facts are
true for the restriction of $u$ to the lower hemisphere. Note that any
(anti-)holomorphic function that is defined on the unit disk and continuous
on the closure is determined by its values on the boundary. Thus, the
function $u_{1}$ is symmetric with respect to the $xy$-plane, as is $u_{2}$.
We conclude that the smooth sections $u=\left( 
\begin{array}{l}
u_{1} \\ 
u_{2}%
\end{array}%
\right) $ in $\ker D$ are exactly functions of $x$ and $y$ alone such that $%
u_{1}$ is holomorphic and $u_{2}$ is antiholomorphic.

We now find the kernel of $D$ restricted to the representation classes of
the $S^{1}$ action. For $k\in \mathbb{Z}_{\geq 0}$ , $w=x+iy$, and $\alpha
\in S^{1}$, 
\begin{eqnarray*}
\psi _{\alpha }\left( 
\begin{array}{c}
w^{k} \\ 
0%
\end{array}%
\right) &=&\psi _{\alpha }\left( 
\begin{array}{c}
\left( \sin \phi \right) ^{k}e^{ik\theta } \\ 
0%
\end{array}%
\right) =\left( 
\begin{array}{c}
\left( \sin \phi \right) ^{k}e^{i\alpha }e^{ik\left( \theta +2\alpha \right)
} \\ 
0%
\end{array}%
\right) =e^{i\left( 2k+1\right) \alpha }\left( 
\begin{array}{c}
w^{k} \\ 
0%
\end{array}%
\right) ; \\
\text{similarly, }\psi _{\alpha }\left( 
\begin{array}{c}
0 \\ 
\overline{w}^{k}%
\end{array}%
\right) &=&e^{-i\left( 2k+1\right) \alpha }\left( 
\begin{array}{c}
0 \\ 
\overline{w}^{k}%
\end{array}%
\right) .
\end{eqnarray*}%
Thus, $\ker D$ is the direct sum of the irreducible representations of $%
S^{1} $ on $\ker D$ corresponding to $\alpha \mapsto $ multiplication by $%
e^{i\left( 2k+1\right) \alpha }$ for $k\in \mathbb{Z}$. Then 
\begin{equation*}
\mathrm{ind}^{\rho _{n}}\left( D\right) =\left\{ 
\begin{array}{ll}
-1 & \text{if }n<0\text{ and }n\text{ is odd} \\ 
1 & \text{if }n>0\text{ and }n\text{ is odd} \\ 
0 & \text{otherwise}%
\end{array}%
\right. ,
\end{equation*}%
where $\rho _{n}$ is the representation $\alpha \mapsto $ multiplication by $%
e^{in\alpha }$.

\subsubsection{The equivariant perturbation $Z$}

Next, we will verify Theorem \ref{transverseGIndexTheorem} by calculating
this same index using an equivariant perturbation $Z$. Let 
\begin{equation*}
Z=\sin \phi \left( 
\begin{array}{cc}
0 & e^{-i\theta } \\ 
e^{i\theta } & 0%
\end{array}%
\right) .
\end{equation*}%
We chose this $Z$ so that at the north and south poles, it will agree with
Clifford multiplication by $\pm i\partial _{\theta }$, as we shall soon see.
The bundle map $Z$ is equivariant; one may check that for any section $u$, $%
Z\left( \psi _{\alpha }u\right) \left( \theta ,\phi \right) =\psi _{\alpha
}\left( Zu\right) \left( \theta ,\phi \right) $. Further, $Z$ is nonsingular
away from the north pole $\phi =0$ and south pole $\phi =\pi $. Next, we see
that $DZ+ZD$ is bounded on sections of the form 
\begin{equation*}
u=\left( 
\begin{array}{c}
f\left( \phi \right) e^{im_{1}\theta } \\ 
g\left( \phi \right) e^{im_{2}\theta }%
\end{array}%
\right) ,
\end{equation*}%
because the coefficient of $\frac{\partial }{\partial \phi }$ in the
expression $DZ+ZD$ is zero. Thus, it is bounded on sections of type $\rho
_{n}$, where $\rho _{n}$ is an irreducible representation of $S^{1}$.

\medskip Since 
\begin{equation*}
D=\left( 
\begin{array}{ll}
0 & -1 \\ 
1 & 0%
\end{array}%
\right) \left( -z\frac{\partial }{\partial x}+x\frac{\partial }{\partial z}%
\right) +\left( 
\begin{array}{ll}
0 & i \\ 
i & 0%
\end{array}%
\right) \left( -z\frac{\partial }{\partial y}+y\frac{\partial }{\partial z}%
\right) ,
\end{equation*}%
at $z=1$ the operator is 
\begin{equation*}
D_{NP}=\left( 
\begin{array}{ll}
0 & 1 \\ 
-1 & 0%
\end{array}%
\right) \frac{\partial }{\partial x}+\left( 
\begin{array}{ll}
0 & -i \\ 
-i & 0%
\end{array}%
\right) \frac{\partial }{\partial y},
\end{equation*}%
and on the whole sphere 
\begin{eqnarray*}
Z &=&i\sin \phi \left( 
\begin{array}{cc}
0 & -ie^{-i\theta } \\ 
-ie^{i\theta } & 0%
\end{array}%
\right) \\
&=&i\left( 
\begin{array}{cc}
0 & -ix-y \\ 
-ix+y & 0%
\end{array}%
\right) =-iy\left( 
\begin{array}{cc}
0 & 1 \\ 
-1 & 0%
\end{array}%
\right) +ix\left( 
\begin{array}{cc}
0 & -i \\ 
-i & 0%
\end{array}%
\right)
\end{eqnarray*}%
At the north pole we define Clifford multiplication as 
\begin{equation*}
c_{NP}\left( a\frac{\partial }{\partial x}+b\frac{\partial }{\partial y}%
\right) :=a\left( 
\begin{array}{ll}
0 & 1 \\ 
-1 & 0%
\end{array}%
\right) +b\left( 
\begin{array}{ll}
0 & -i \\ 
-i & 0%
\end{array}%
\right) ,
\end{equation*}%
then 
\begin{equation*}
Z=ic_{NP}\left( -y\frac{\partial }{\partial x}+x\frac{\partial }{\partial y}%
\right) =ic_{NP}\left( \partial _{\theta }\right) .
\end{equation*}%
Similarly, at the south pole, 
\begin{equation*}
D_{SP}=\left( 
\begin{array}{ll}
0 & -1 \\ 
1 & 0%
\end{array}%
\right) \frac{\partial }{\partial x}+\left( 
\begin{array}{ll}
0 & i \\ 
i & 0%
\end{array}%
\right) \frac{\partial }{\partial y},
\end{equation*}%
and we define Clifford multiplication as 
\begin{equation*}
c_{SP}\left( a\frac{\partial }{\partial x}+b\frac{\partial }{\partial y}%
\right) :=a\left( 
\begin{array}{ll}
0 & -1 \\ 
1 & 0%
\end{array}%
\right) +b\left( 
\begin{array}{ll}
0 & i \\ 
i & 0%
\end{array}%
\right) ,
\end{equation*}%
so that 
\begin{equation*}
Z=ic_{SP}\left( y\frac{\partial }{\partial x}-x\frac{\partial }{\partial y}%
\right) =ic_{SP}\left( -\partial _{\theta }\right) .
\end{equation*}

We now use Theorem \ref{transverseGIndexTheorem} to calculate the index $%
\mathrm{ind}^{\rho _{n}}\left( D\right) $. We consider sections of type $%
\rho _{n}$. At the north pole $z=1$, we have%
\begin{eqnarray*}
\tau _{1} &=&1;~\mathbf{b}=b_{1}=n;~k_{11}=\kappa _{1}=2 \\
a_{11} &=&-1;~a_{21}=1;~\varepsilon _{11}=1;~\varepsilon _{21}=-1
\end{eqnarray*}%
Then the contribution to the index at the north pole is%
\begin{eqnarray*}
\mathrm{sign}\left( 1\right) k_{1}\left( NP\right) +\mathrm{sign}\left(
2\right) k_{2}\left( NP\right) &=&0-k_{2}\left( NP\right) \\
&=&-\#\left\{ \left. m\in \mathbb{Z\,~}\right\vert \text{ }m\leq 0\text{, }%
\varepsilon _{21}=-1\text{, and }mk_{11}=\,a_{21}+b_{1}\right\} \\
&=&-\#\left\{ \left. m\in \mathbb{Z\,~}\right\vert \text{ }m\leq 0\text{, }%
\varepsilon _{21}=-1\text{, and }2m=\,1+n\right\} \\
&=&-1\text{ if }n<0\text{ and }n\text{ is odd.}
\end{eqnarray*}%
At the south pole, with the orientation reversed on both the surface and on
the group $S^{1}$,%
\begin{eqnarray*}
\tau _{1} &=&1;~\mathbf{b}=b_{1}=-n;~k_{11}=\kappa _{1}=2 \\
a_{11} &=&1;~a_{21}=-1;~\varepsilon _{11}=-1;~\varepsilon _{21}=1
\end{eqnarray*}%
Then the contribution to the index at the south pole is%
\begin{eqnarray*}
\mathrm{sign}\left( 1\right) k_{1}\left( SP\right) +\mathrm{sign}\left(
2\right) k_{2}\left( SP\right) &=&k_{1}\left( NP\right) -0 \\
&=&\#\left\{ \left. m\in \mathbb{Z\,~}\right\vert \text{ }m\leq 0\text{, }%
\varepsilon _{11}=-1\text{, and }mk_{11}=\,a_{11}+b_{1}\right\} \\
&=&\#\left\{ \left. m\in \mathbb{Z\,~}\right\vert \text{ }m\leq 0\text{, }%
\varepsilon _{11}=-1\text{, and }2m=\,1-n\right\} \\
&=&1\text{ if }n>0\text{ and }n\text{ is odd.}
\end{eqnarray*}%
Then, as expected, we have%
\begin{equation*}
\mathrm{ind}^{\rho _{n}}\left( D\right) =\left\{ 
\begin{array}{ll}
-1 & \text{if }n<0\text{ and }n\text{ is odd} \\ 
1 & \text{if }n>0\text{ and }n\text{ is odd} \\ 
0 & \text{otherwise}%
\end{array}%
\right. .
\end{equation*}

\section{Appendix}

\subsection{Proof of Lemma \protect\ref{DZ+ZD Lemma}}

\begin{proof}
We have 
\begin{equation*}
\left( D_{s}\right) ^{2}-D^{2}=s\left( ic\left( V\right) \circ D+iD\circ
c\left( V\right) \right) +s^{2}\left| V\right| ^{2}.
\end{equation*}
Write $V=\sum V_{j}e_{j}$ in terms of a local orthonormal frame $%
e_{1},e_{2},...$ of the tangent bundle, corresponding to geodesic normal
coordinate vector fields $e_{j}=\partial _{j}$ at the origin of the
coordinate system. At the origin of the coordinate system, we have 
\begin{eqnarray*}
ZD+DZ &=&ic\left( V\right) \circ D+iD\circ c\left( V\right) \\
&=&i\sum_{j,k}V_{j}c\left( e_{j}\right) c\left( e_{k}\right) \nabla
_{k}+i\sum_{j,k}c\left( e_{k}\right) \nabla _{k}\circ V_{j}c\left(
e_{j}\right) .
\end{eqnarray*}
Since $\nabla _{k}e_{j}=0$ at the origin, 
\begin{eqnarray*}
ZD+DZ &=&i\sum_{j,k}V_{j}c\left( e_{j}\right) c\left( e_{k}\right) \nabla
_{k}+i\sum_{j,k}V_{j}c\left( e_{k}\right) c\left( e_{j}\right) \nabla
_{k}+i\sum_{j,k}c\left( e_{k}\right) c\left( e_{j}\right) \partial _{k}V_{j}
\\
&=&-2i\sum_{j}V_{j}\nabla _{j}-i\sum_{j}\partial _{j}V_{j}+i\sum_{j\neq
k}c\left( e_{k}\right) c\left( e_{j}\right) \partial _{k}V_{j}\text{, }
\end{eqnarray*}
since $c\left( e_{j}\right) c\left( e_{k}\right) +c\left( e_{k}\right)
c\left( e_{j}\right) =-2\delta _{jk}$. Then 
\begin{gather}
ZD+DZ=-2i\nabla _{V}-i\left( \mathrm{div}\left( V\right)
-\sum_{j}V_{j}\sum_{k\neq j}\left\langle \nabla _{k}e_{j},e_{k}\right\rangle
\right) +i\sum_{j\neq k}c\left( e_{k}\right) c\left( e_{j}\right) \partial
_{k}V_{j}  \notag \\
=-2i\nabla _{V}-i\mathrm{div}\left( V\right) +i\sum_{j\neq k}c\left(
e_{k}\right) c\left( e_{j}\right) \partial _{k}V_{j}\text{ since }\nabla
_{k}e_{j}=0\text{ at the origin}  \notag \\
=-2i\nabla _{V}-i\mathrm{div}\left( V\right) +ic\left( d\left( V^{*}\right)
\right) ,  \notag
\end{gather}
where by $c\left( d\left( V^{*}\right) \right) $ we imply that we have used
the inverse of the symbol map $\sigma $ to convert the two-form $d\left(
V^{*}\right) $ to a Clifford algebra element. For example, $\sigma \left(
e_{1}e_{2}\right) =c\left( e_{1}\right) c\left( e_{2}\right) \mathrm{1}%
=\left( dx_{1}\wedge -dx_{1}\lrcorner \right) \left( dx_{2}\wedge
-dx_{2}\lrcorner \right) \mathrm{1}=dx_{1}\wedge dx_{2}$, so we define $%
c\left( dx_{1}\wedge dx_{2}\right) =c\left( e_{1}e_{2}\right) $ at the
origin. Now, since the last expression is coordinate-free, we conclude that 
\begin{equation*}
ZD+DZ=-2i\nabla _{V}-i\mathrm{div}\left( V\right) +ic\left( d\left(
V^{*}\right) \right)
\end{equation*}
at all points.
\end{proof}

\subsection{Proof of the statement in Example \protect\ref{SpinorDZ+ZD}}

\begin{proof}
A Killing vector field $X$ can be lifted to a vector field $\overline{X}$ on
the frame bundle, so that $\overline{X}$ covers $X$ and is invariant under
the $SO\left( n\right) $ bundle. The vector field $\overline{X}$ lifts
uniquely to a vector field $\widehat{X}$ on the principal spin bundle $%
\widetilde{F}$. Thus it acts on any bundle associated to $\widetilde{F}$,
such as the spin bundle. Let $X$ be an infinitesimal isometry. If $g$ is the
metric tensor, then $\mathcal{L}_{X}g=0$. If $Y$, $Z$ are any two tensor
fields of the same type, then 
\begin{eqnarray*}
X\left\langle Y,Z\right\rangle &=&\left\langle \mathcal{L}
_{X}Y,Z\right\rangle +\left\langle Y,\mathcal{L}_{X}Z\right\rangle \\
&=&\left\langle \nabla _{X}Y,Z\right\rangle +\left\langle Y,\nabla
_{X}Z\right\rangle .
\end{eqnarray*}
Thus $A_{X}=\mathcal{L}_{X}-\nabla _{X}$ is skew-symmetric and of degree
zero, since $\left\langle A_{X}Y,Z\right\rangle =-\left\langle
Y,A_{X}Z\right\rangle .$ Hence its action on $\Gamma \left( TM\right) $
comes from the endomorphism (also called $A_{X}$) of $TM$. Choose a basis $%
\left\{ e_{i}\right\} $ of $T_{x}M$, and we may identify $A_{X}$ with the
element $a_{X}\in \mathfrak{o}\left( n\right) $ by identifying $T_{x}M$ with 
$\mathbb{R}^{n}$ using the basis . Under this identification the
antisymmetric matrix $a_{X}=\left( \left( a_{X}\right) _{ij}\right) $
corresponds to an endomorphism $A_{X}=$ $\frac{1}{4}\sum \left( a_{X}\right)
_{ij}e_{i}e_{j}$ (see Lemma 4.8 in \cite{Roe} for calculations).

Let $\lambda :\mathrm{Spin}\left( n\right) \rightarrow SO\left( n\right) $
be the double cover, and let $d\lambda :\mathfrak{spin}\left( n\right)
\rightarrow \mathfrak{o}\left( n\right) $ be the differential map on the Lie
algebras. Observe that $\mathfrak{spin}\left( n\right) \cong \mathrm{Cl}
_{2}\left( \mathbb{R}^{n}\right) $, and the Lie bracket induced on $\mathrm{%
\ \ \ Cl}_{2}\left( \mathbb{R}^{n}\right) $ is $\left[ a,b\right] =ab-ba$
(using Clifford multiplication). For all $v\in \mathbb{R}^{n}$ , $z\in 
\mathrm{Spin}\left( n\right) $, 
\begin{equation*}
d\lambda \left( z\right) \left( v\right) =zv-vz,
\end{equation*}
where $z$ is thought of as an element of $\mathrm{Cl}_{2}\left( \mathbb{R}%
^{n}\right) $ and $v$ is thought of as an element of $\mathrm{Cl}_{1}\left( 
\mathbb{R}^{n}\right) $. Hence 
\begin{eqnarray*}
d\lambda \left( \frac{1}{4}\sum \left( a_{X}\right) _{ij}e_{i}e_{j}\right)
\left( v\right) &=&\left[ \frac{1}{4}\sum \left( a_{X}\right)
_{ij}e_{i}e_{j},v\right] \\
&=&a_{X}v,
\end{eqnarray*}
so that 
\begin{equation*}
A_{X}=\mathcal{L}_{X}-\nabla _{X}=\frac{1}{4}\sum \left( a_{X}\right)
_{ij}e_{i}e_{j}.
\end{equation*}
Next, given a Killing field $X$ and vector field $Y$, 
\begin{eqnarray*}
A_{X}Y &=&\mathcal{L}_{X}Y-\nabla _{X}Y \\
&=&\left[ X,Y\right] -\left[ X,Y\right] -\nabla _{Y}X \\
&=&-\nabla _{Y}X.
\end{eqnarray*}
Thus, given any vector field $Z$, 
\begin{equation*}
\left\langle A_{X}Y,Z\right\rangle =-\left\langle \nabla
_{Y}X,Z\right\rangle ,
\end{equation*}
so $A_{X}=-\left( \nabla X\right) ^{\#}$. This implies 
\begin{equation*}
\left( a_{X}\right) _{ij}=-\left\langle \nabla _{e_{i}}X,e_{j}\right\rangle .
\end{equation*}
Thus, 
\begin{eqnarray*}
A_{X} &=&\mathcal{L}_{X}-\nabla _{X}=-\frac{1}{4}\sum \left\langle \nabla
_{e_{i}}X,e_{j}\right\rangle e_{i}e_{j} \\
&=&-\frac{1}{4}\sum e_{i}\left( \left\langle \nabla
_{e_{i}}X,e_{j}\right\rangle e_{j}\right) \\
&=&-\frac{1}{4}\sum e_{i}\left( \nabla _{e_{i}}X\right) \\
&=&-\frac{1}{4}\sum e_{i}\left( \partial _{i}X_{j}\right) e_{j}\text{ if }%
\left\{ e_{j}\right\} \text{ is isochronous} \\
&=&-\frac{1}{4}c\left( d\left( X^{*}\right) \right)
\end{eqnarray*}
We have therefore that 
\begin{equation*}
\mathcal{L}_{X}=\nabla _{X}-\frac{1}{4}c\left( d\left( X^{*}\right) \right)
\end{equation*}
if $X$ is a Killing vector field.
\end{proof}


\begin{thebibliography}{99}
\bibitem{A} M. F. Atiyah, \textit{Elliptic operators and compact groups},
Lecture Notes in Math. \textbf{401}, Springer-Verlag, Berlin, 1974.

\bibitem{AB1} M. F. Atiyah and R. Bott, \emph{A Lefschetz fixed point
formula for elliptic complexes I}, Ann. of Math. (2) \textbf{86}(1967),
374--407.

\bibitem{AB2} M. F. Atiyah and R. Bott, \emph{A Lefschetz fixed point
formula for elliptic complexes II}, Ann. of Math. (2) \textbf{88}(1968),
451--491.

\bibitem{AH} M. F. Atiyah and F. Hirzebruch, \emph{Spin-manifolds and group
actions}, in \textit{Essays on Topology and Related Topics (M\'{e}moires d%
\'{e}di\'{e}s \`{a} Georges de Rham)}, Springer-Verlag, New York, 1970,
18--28.

\bibitem{ASe} M. F. Atiyah and G. B. Segal, \emph{The index of elliptic
operators: II}, Ann. of Math. (2) \textbf{87}(1968), 531--545.

\bibitem{B-G-V} N. Berline, E. Getzler, and M. Vergne, \textit{Heat kernels
and Dirac operators}, Grundlehren der Mathematischen Wissenschaften \textbf{%
\ \ \ \ 298}, Springer-Verlag, Berlin, 1992.

\bibitem{Be-V1} N. Berline and M. Vergne, \emph{The Chern character of a
transversally elliptic symbol and the equivariant index}, Invent. Math. 
\textbf{124}(1996), no. 1-3, 11-49.

\bibitem{Be-V2} N. Berline and M. Vergne, \emph{L'indice \'{e}quivariant des
op\'{e}rateurs transversalement elliptiques}, Invent. Math. \textbf{124}%
(1996), no. 1-3, 51-101.

\bibitem{Bra} M. Braverman, \emph{Index theorem for equivariant Dirac
operators on non-compact manifolds}, K-Theory \textbf{27}(2002), 61-101.

\bibitem{BrotD} T. Br\"{o}cker and T. tom Dieck, \emph{Representations of
Compact Lie Groups}, Graduate Texts in Math., Springer-Verlag, New York,
1985.

\bibitem{BruKRich} J. Br\"{u}ning, F.\ W. Kamber, and K. Richardson, \emph{%
The equivariant index of transversally elliptic operators}, in preparation.

\bibitem{G-S} V. Guillemin and S. Sternberg, \emph{Geometric quantization
and multiplicities of group representations}, Invent. Math. \textbf{67}%
(1982), 515--538.

\bibitem{Kawas2} T. Kawasaki, \emph{The index of elliptic operators over }$V$%
\emph{-manifolds}, Nagoya Math. J. \textbf{84} (1981), 135--157.

\bibitem{Par} P. E. Paradan, \emph{Localization of the Riemann-Roch character%
}, J. Funct. Anal. \textbf{187}(2001), 442-509.

\bibitem{Pro-Rich} I. Prokhorenkov, K. Richardson, \emph{Perturbations of
Dirac operators}, J. Geom. Phys. \textbf{57}(2006), 297-321.

\bibitem{Roe} J. Roe,\textit{\ Elliptic operators, topology, and asymptotic
methods}, Pitman Research Notes in Math. \textbf{179}, Longman Scientific
and Technical, Harlow, 1988.

\bibitem{TZ} Y. Tian and W. Zhang, \emph{An analytic proof of the geometric
quantization conjecture of Guillemin-Sternberg}, Invent. Math. \textbf{132}%
(1998), no. 2, 229--259.

\bibitem{Shu1} M. A. Shubin, \textit{Semiclassical asymptotics on covering
manifolds and Morse inequalities}, Geom. Funct. Anal.\textbf{\ 6} (1996),
no. 2, 370--409.

\bibitem{Wit1} E. Witten, \textit{Supersymmetry and Morse Theory}, J.
Differ. Geometry, \textbf{17}, 661-692 (1982).

\bibitem{Wit3} E. Witten, \textit{Index of Dirac operators}, Quantum fields
and strings: a course for mathematicians, Vol. 1, 475--511, Amer. Math.
Soc., Providence, RI, 1999.
\end{thebibliography}
\end{document}